\documentclass[reqno,letterpaper]{amsart}
\usepackage{import}
\usepackage{amsmath,amsthm,amssymb}
\usepackage{thmtools}
\usepackage{url}
\usepackage{comment}
\usepackage{graphicx}
\usepackage{subcaption}
\usepackage[numbers,sort&compress]{natbib}
\RequirePackage[usenames,dvipsnames]{color}
\usepackage[colorlinks=true,linkcolor=blue,citecolor=blue,pdfborder={0 0 0}]{hyperref}
\usepackage[capitalize]{cleveref}
\usepackage{datetime}

\usepackage[arxiv]{optional}

\graphicspath{ {../figures/}{figures/} }

\usepackage{jhh-misc2}

\declaretheorem[style=theorem,name=Assumption]{assumption}
\renewcommand*{\theassumption}{\arabic{section}.\Alph{assumption}}

\declaretheorem[style=remark,qed=$\square$,name=Remark,sibling=nthm]{remark}
\declaretheorem[numbered=no,style=definition,qed=$\square$,name=Example]{example}

\crefname{nlem}{Lemma}{Lemmas}
\crefname{nprop}{Proposition}{Propositions}
\crefname{ncor}{Corollary}{Corollaries}
\crefname{nthm}{Theorem}{Theorems}
\crefname{assumption}{Assumption}{Assumptions}
\crefname{property}{Property}{Properties}
\crefname{remark}{Remark}{Remarks}

\numberwithin{equation}{section}

\allowdisplaybreaks[3]

\newcommand{\gen}[1]{\mcA_{#1}}
\newcommand{\ball}{\mathbb{B}}
\newcommand{\id}{\textrm{id}}
\newcommand{\distDiff}{\distNamed{Diff}}

\def\figwidth{.45}

\usepackage{autonum}

\begin{document}

\title[Approximate Diffusions]{Quantifying the accuracy of approximate \\ diffusions and Markov chains}

\author[J.~H.~Huggins]{Jonathan H.~Huggins}
\address{Computer Science and Artificial Intelligence Laboratory (CSAIL) \\ Massachusetts Institute of Technology}
\urladdr{http://www.jhhuggins.org/}
\email{jhuggins@mit.edu}

\author[J.~Zou]{James Zou}
\address{Stanford University}
\urladdr{http://sites.google.com/site/jamesyzou/}
\email{jamesyzou@gmail.com}

\date{\today}

\begin{abstract}
Markov chains and diffusion processes are indispensable tools in machine learning 
and statistics that are used for inference, sampling, and modeling. 
With the growth of large-scale datasets, the computational cost associated with 
simulating these stochastic processes can be considerable, and many algorithms have been 
proposed to approximate the underlying Markov chain or diffusion. 
A fundamental question is how the computational savings trade off against the 
statistical error incurred due to approximations. 
This paper develops general results that address this question. 
We bound the Wasserstein distance between the equilibrium distributions of two 
diffusions as a function of their mixing rates and the deviation in their drifts. 
We show that this error bound is tight in simple Gaussian settings. 
Our general result on continuous diffusions can be discretized to provide 
insights into the computational--statistical trade-off of Markov chains. 
As an illustration, we apply our framework to derive finite-sample error 
bounds of approximate unadjusted Langevin dynamics. 
We characterize computation-constrained settings where, by using 
fast-to-compute approximate gradients in the Langevin dynamics, we obtain 
more accurate samples compared to using the exact gradients. 
Finally, as an additional application of our approach, we quantify 
the accuracy of approximate zig-zag sampling. 
Our theoretical analyses are supported by simulation experiments.

\end{abstract}

\maketitle

\newcommand{\suppmat}{\opt{XXXXX}{Appendix}\opt{arxiv}{Appendix}}

\section{Introduction}

Markov chains and their continuous-time counterpart, diffusion processes, are ubiquitous 
in machine learning and statistics, forming a core component of the inference and 
modeling toolkit. 
Since faster convergence enables more efficient sampling and inference, a large and 
fruitful literature has investigated how quickly these stochastic processes converge 
to equilibrium.  
However, the tremendous growth of large-scale machine learning datasets -- in areas such as social 
network analysis, vision, natural language processing and bioinformatics -- have created
new inferential challenges.
The large-data setting highlights the need for stochastic processes that 
are not only accurate (as measured by fast convergence to the target distribution),
but also computationally efficient to simulate.
These computational considerations have led to substantial research efforts into 
approximating the underlying stochastic processes with new processes that are more
computationally efficient \citep{Welling:2011,Bardenet:2015,Ge:2015}.

As an example, consider using Markov chain Monte Carlo (MCMC) to sample from
a posterior distribution. 
In standard algorithms, each step of the Markov chain involves calculating a 
statistic that depends on all of the observed data (e.g.~a likelihood ratio 
to set the rejection rate in Metropolis-Hastings or a gradient of the 
log-likelihood as in Langevin dynamics). 
As data sets grow larger, such calculations increasingly become the computational bottleneck. 
The need for more scalable sampling algorithms has led to the development of 
Markov chains which only approximate the desired statistics at each 
step -- for example, by approximating the gradient or sub-sampling the data -- and hence are computationally more 
efficient~\citep{Welling:2011,Chen:2014,Maclaurin:2014,Korattikara:2014,Bardenet:2014,Bardenet:2015,Quiroz:2015,Ge:2015}. 
The trade-off is that the approximate chain often does not converge to the desired 
equilibrium distribution, which, in many applications, could be the posterior distribution 
of some latent parameters given all of the observed data. 
Therefore, a central question of both theoretical and practical importance is how to quantify 
the deviation between the equilibrium distribution that the approximate chain converges to 
and the desired distribution targeted by the original chain. 
Moreover, we would like to understand, given a fixed computational budget, how to design 
approximate chains that generate the most accurate samples.

\textbf{Our contributions.} 
In this paper, we develop general results to quantify the accuracy of approximate diffusions 
and Markov chains and apply these results to characterize the computational--statistical 
trade-off in specific algorithms. 
Our starting point is continuous-time diffusion processes because these are the  
objects which are discretized to construct many sampling algorithms, such as the 
unadjusted and Metropolis-adjusted Langevin algorithms~\citep{Roberts:1996} and Hamiltonian 
Monte Carlo~\citep{Neal:2011}. 
Given two diffusion processes, we bound the deviation in their equilibrium distributions 
in terms of the deviation in their drifts and the rate at which the diffusion 
mixes (\cref{thm:approx-grad-distance}). 
Moreover, we show that this bound is tight for certain Gaussian target distributions.
These characterizations of diffusions are novel
and are likely of more general interest beyond the inferential settings we consider. 
We apply our general results to derive a finite-sample error bound on a specific 
unadjusted Langevin dynamics algorithm (\cref{thm:ula-accuracy-main}). 
Under computational constraint, the relevant trade-off here is between computing 
the exact log-likelihood gradient for few iterations or computing an approximate 
gradient for more iterations. 
We characterize settings where the approximate Langevin dynamics produce more
accurate samples from the true posterior. 
We illustrate our analyses with simulation results. 
In addition, we apply our approach to quantify the accuracy of 
approximations to the zig-zag process, a recently-developed non-reversible sampling scheme.

\textbf{Paper outline.} 
We introduce the basics of diffusion processes and other preliminaries in \cref{sec:preliminaries}. 
\cref{sec:main} discusses the main results on bounding the error between an exact and perturbed diffusion. 
We describe the main ideas behind our analyses in \cref{sec:overview}; 
all the detailed proofs are deferred to the \suppmat. 
\cref{sec:trade-off} applies the main results to derive finite sample error 
bounds for unadjusted Langevin dynamics and illustrates the 
computational--statistical trade-off. 
\cref{sec:pdmps} extends our main results to quantify the accuracy of
approximate piecewise deterministic Markov processes, including the zig-zag process.
Numerical experiments to complement the theory are provided in \cref{sec:experiments}. 
We conclude with a discussion of how our results connect to the relevant 
literature and suggest directions for further research. 

\section{Diffusions and preliminaries}
\label{sec:preliminaries}

Let $\mcX = \reals^d$ be the parameter space and let $\pi$ be a probability density 
over $\reals^d$ (e.g.~it can be the posterior distribution of some latent parameters given data). 
A Langevin diffusion is characterized by the stochastic differential equation
\[
\dee X_{t} = \grad \log \pi(X_t)\,\dee t + \sqrt{2} \dee W_{t}, \label{eqn:LD}
\]
where $X_t\in \reals^d$ and $W_t$ is a standard Brownian motion. 
The intuition is that $X_t$ undergoes a biased random walk in which it is more likely to 
move in directions that increase the density. 
Under appropriate regularity conditions, as $t \rightarrow \infty$, the distribution of $X_t$ converges to $\pi$. 
Thus, simulating the Langevin diffusion provides a powerful framework to sample 
from the target $\pi$.
To implement such a simulation, we need to discretize the continuous diffusion into
finite-width time steps. 
For our main results, we focus on analyzing properties of the underlying diffusion processes. 
This allows us to obtain general results which are independent of any particular discretization scheme. 

Beyond Langevin dynamics, more general diffusions can take the form
\[
\dee X_{t} = b(X_t)\,\dee t + \sqrt{2} \dee W_{t}, \label{eqn:diffusion}
\]
where $b: \reals^d \rightarrow \reals^d$ is the drift and is not necessarily the 
gradient of some log-density.\footnote{All of our results can be extended to more 
general diffusions on a domain $\mathcal{X} \subseteq \reals^d$, 
$\dee X_{t} = b(X_t) + \Sigma\,\dee W_{t} - n_t L(\dee t)$, where $\Sigma$ is the 
covariance of the Brownian motion, and $n_t L$ captures the reflection forces 
at the boundary $\partial \mathcal{X}$. 
To keep the exposition simple, we focus on the simpler diffusion in the main text.}
Furthermore, we can analyze other continuous-time Markov processes such as
piecewise deterministic Markov processes (PDMPs).
For example, Hamiltonian Monte Carlo~\citep{Neal:2011} can be viewed
as approximating a PDMP and the zig-zag process is a recently-developed non-reversible
PDMP designed for large Bayesian inference (see~\cref{sec:pdmps}).

In many large-data settings, computing the drift $b(X_t)$ in \cref{eqn:diffusion} can 
be expensive; for example, computing $b(X_{t}) = \grad \log \pi(X_t)$ requires using
all of the data and may involve evaluating a complex function such as a differential 
equation solver. 
Many recent algorithms have been proposed where we replace $b$ with an approximation $\tb$. 
Such an approximation changes the underlying diffusion process to 
\[
\dee \tX_{t} = \tb(\tX_t)\,\dee t + \sqrt{2} \dee \tW_{t}, \label{eqn:diffusion-approx}
\]
where $\tW_{t}$ is a standard Brownian motion.
In order to understand the quality of different approximations, we need to quantify how the 
equilibrium distribution of \cref{eqn:diffusion} differs from the equilibrium distribution of \cref{eqn:diffusion-approx}. 
We use the standard \emph{Wasserstein metric} to measure this distance. 

\begin{defn}
The \emph{Wasserstein distance} between distributions $\pi$ and $\tpi$ is 
\[
d_{\mcW}(\pi, \tpi) = \sup_{\phi \in C_{L}(\reals^{d})} |E_{\pi}[\phi] - E_{\tpi}[\phi]|,
\]
where $C_{L}(\reals^{d})$ is the set of continuous functions $\phi : \reals^d \to \reals$ with Lipschitz 
constant $\|\phi\|_L \leq 1$.\footnote{Recall that the Lipschitz constant of function
$\phi :\reals^d \to \reals$ is $\|\phi\|_L \defined  \sup_{x, y \in \reals^{d}} \frac{\|\phi(x) - \phi(y)\|_2}{\|x-y\|_2}$.}
\end{defn}

The distance between $\pi$ and $\tpi$ should depend on how good the drift approximation is,
which can be quantified by $\|b - \tb\|_2$.\footnote{For a function $\phi : \reals^{n} \to \reals^{m}$, define $\|\phi\|_{2} \defined \sup_{x \in \reals^{n}}\|\phi(x)\|_{2}$.} 
It is also natural for the distance to depend on how quickly the original diffusion with drift $b$ mixes,
since the faster it mixes, the less time there is for the error to accumulate. 
Geometric contractivity is a useful property which quantifies fast-mixing diffusions. 
For each $x \in \reals^{d}$, let $\mu_{x,t}$ denote the law of $X_{t} \given X_{0} = x$. 

\begin{assumption}[\textbf{Geometric contractivity}] \label{asm:exponential-contractivity}
There exist constants $C > 0$ and $0 < \rho < 1$ such that
for all $x, x' \in \reals^{d}$, 
\[
d_{\mcW}(\mu_{x,t}, \mu_{x',t}) \le C\|x - x'\|_{2}\rho^{t}.
\]
\end{assumption}

Geometric contractivity holds in many natural settings. 
Recall that a twice continuously-differentiable function $\phi$ is \emph{$k$-strongly concave} if for all $x, x' \in \mathbb{R}^d$
\[
(\grad \phi(x) - \grad \phi(x')) \cdot (x - x') \le -k \|x - x'\|_{2}^{2}.
\label{eqn:exponential-erg}
\]
When $b = \grad \log \pi$ and $\log \pi$ is $k$-strongly concave, the diffusion 
is exponentially ergodic with $C = 1$ and $\rho = e^{-k}$ (this can be shown using 
standard coupling arguments~\citep{Bolley:2012}). 
In fact, exponential contractivity also follows if \cref{eqn:exponential-erg} is satisfied 
when $x$ and $x'$ are far apart and $\log \pi$ has ``bounded convexity'' 
when $x$ and $x'$ are close together~\citep{Eberle:2015}.
Alternatively, \citet{Hairer:2009} provides a Lyapunov function-based approach to proving exponential contractivity.

To ensure that the diffusion and the approximate diffusion are well-behaved, 
we impose some standard regularity properties. 

\begin{assumption}[\textbf{Regularity conditions}] \label{asm:regularity}
Let $\pi$ and $\tpi$ denote the stationary distributions of the 
diffusions in \cref{eqn:diffusion} and \cref{eqn:diffusion-approx}, respectively. 
\benum
\item The target density satisfies $\pi \in C^{2}(\reals^d, \reals)$ and $\int x^{2}\pi(\dee x) < \infty$.
The drift satisfies $b \in C^{1}(\reals^d, \reals^{d})$ and $\|b\|_{L} < \infty$. 
\item The approximate drift satisfies $\tb \in C^{1}(\reals^d, \reals^{d})$ and $\|\tb\|_{L} < \infty$. 
\item If a function $\phi \in C(\reals^d, \reals)$ is $\pi$-integrable then it is $\tpi$-integrable.
\eenum
\end{assumption}

Here $C^k(\reals^m, \reals^n)$ denotes the set of $k$-times continuously differentiable functions from 
$\reals^m$ to $\reals^n$ and $C(\reals^m, \reals^n)$ is the set of all Lebesgue-measurable function 
from $\reals^m$ to $\reals^n$. 
The only notable regularity condition is (3). 
In the \suppmat, we discuss how to verify it and why it can safely be treated
as a mild technical condition. 

\section{Main results}\label{sec:main}

We can now state our main result, which quantifies the deviation in the equilibrium distributions
of the two diffusions in terms of the mixing rate and the difference between the diffusions' drifts. 

\bnthm[Error induced by approximate drift] \label{thm:approx-grad-distance}
Let $\pi$ and $\tpi$ denote the invariant distributions of the 
diffusions in \cref{eqn:diffusion} and \cref{eqn:diffusion-approx}, respectively. 
If the diffusion \cref{eqn:diffusion} is exponentially ergodic with parameters $C$ and $\rho$, 
the regularity conditions of \cref{asm:regularity} hold, and $\|b - \tb\|_{2} \le \eps$, then
\[
d_{\mcW}(\pi, \tpi) \le \frac{C\eps}{\log(1/\rho)}. \label{eq:det-gradient-approx-bound}
\]
\enthm

\begin{remark}[\emph{Coherency of the error bound}]
To check that the error bound of \cref{eq:det-gradient-approx-bound} has coherent dependence 
on its parameters, consider the following thought experiment. 
Suppose we change the time scale of the diffusion from $t$ to $s = at$ for some $a > 0$. 
We are simply \emph{speeding up} or \emph{slowing down} the diffusion process depending 
on whether $a>1$ or $a<1$. 
Changing the time scale does not affect the equilibrium distribution and hence 
$d_{\mcW}(\pi, \tpi)$ remains unchanged. 
After time $s$ has passed, the exponential contraction is $\rho^{at}$ and hence the 
effective contraction constant is $\rho^a$ instead of $\rho$. 
Moreover, the drift at each location is also scaled by $a$ and hence the drift 
error is $\epsilon a$. 
The scaling $a$ thus cancels out in the error bound, which is desirable since the 
error should be independent of how we set the time scale.   
\end{remark} 

\begin{remark}[\emph{Tightness of the error bound}]
We can choose $b$ and $\tb$ such that the bound in \cref{eq:det-gradient-approx-bound} 
is an equality, thus showing that, under the assumptions considered, 
\cref{thm:approx-grad-distance} cannot be improved.  
Let $\pi(x) = \distNorm(x; \mu, \sigma^{2}I)$ be the Gaussian density with mean 
$\mu \in \reals^{d}$ and covariance matrix $\sigma^{2}I$ and let 
$\tpi(x) = \distNorm(x; \tmu, \sigma^{2}I)$. 
The Wasserstein distance between two Gaussians with the same covariance is the 
distance between their means, so $d_{\mcW}(\pi, \tpi) = \|\mu - \tmu\|_{2}$.
Consider the corresponding diffusions where  $b = \grad \log \pi$ and  $\tb = \grad \log \tpi$. 
We have that for any $x \in \reals^{d}$, $\|b(x) - \tb(x)\|_{2} = \sigma^{-2}\|\mu - \tmu\|_{2} =: \eps$. 
Furthermore, the Hessian is $H[\log \pi] = -\sigma^{-2}I$, which implies that $b$ is 
$\sigma^{-2}$-strongly concave. 
Therefore, per the discussion in \cref{sec:preliminaries}, 
exponential contractivity holds with $C = 1$ and $\rho = e^{-\sigma^{-2}}$.
We thus conclude that 
\[
\frac{C \eps}{\log(1/\rho)} 
= \frac{\sigma^{-2}\|\mu - \tmu\|_{2}}{\sigma^{-2}} 
= \|\mu - \tmu\|_{2} = d_{\mcW}(\pi, \tpi). 
\]
and hence the bound of \cref{thm:approx-grad-distance} is tight in this setting. 
\end{remark}

\cref{thm:approx-grad-distance} assumes that the approximate drift is a deterministic 
function and that the error in the drift is uniformly bounded.
We can generalize the results of \cref{thm:approx-grad-distance} to allow for
the approximate diffusion to use stochastic drift with non-uniform drift error.
We will see that only the expected magnitude of the
drift bias affects the final error bound. 
Let $\tb(\tX_{t}, \tY_{t})$ denote the approximate drift, which is now a function of 
both the current location $\tX_{t}$ and an independent diffusion $\tY_{t} \in \reals^{\ell}$:
\[
&\dee \tX_{t} = (\tb(\tX_{t}, \tY_{t}))\,\dee t + \sqrt{2}\dee \tW_{t}^{X} \label{eq:stochastic-gradient-diffusion} \\
&\dee \tY_{t} = b_{aux}(\tY_{t})\,\dee t + \Sigma\,\dee\tW^{Y}_{t},
\]
where $\Sigma$ is an $\ell \times \ell$ matrix and 
the notation $\tW_{t}^{X}$ and $\tW_{t}^{Y}$ highlights that the Brownian motions 
in $\tX_t$ and $\tY_t$ are independent. 
Let $\tpi_{Z}$ denote the stationary distribution of $\tZ_{t} \defined (\tX_{t}, \tY_{t})$.
For measure $\mu$ and function $f$, we write $\mu(f) \defined \int f(x)\mu(\dee x)$ to reduce clutter.  
We can now state a generalization of \cref{thm:approx-grad-distance}. 

\bnthm[Error induced by stochastic approximate drift] \label{thm:stochastic-approx-grad-distance}
Let $\pi$ and $\tpi$ denote the invariant distributions of the 
diffusions in \cref{eqn:diffusion,eq:stochastic-gradient-diffusion}, respectively.
Assume that there exists a measurable function $\eps \in C(\reals^d, \reals_{+})$
such that for $(\tX, \tY) \dist \tpi_{Z}$ and for all $x \in \reals^d$,
\[
\|b(x) - \EE[\tb(\tX, \tY) \given \tX = x]\|_{2} \le \eps(x).
\]
If the diffusion \cref{eqn:diffusion} is exponentially ergodic and the regularity conditions
of \cref{asm:regularity} hold, then
\[
d_{\mcW}(\pi, \tpi) \le \frac{C\,\tpi(\eps)}{\log(1/\rho)}.
\]
\enthm

Whereas the bound of \cref{thm:approx-grad-distance} is proportional to 
the deterministic drift error $\epsilon$, the bound for the diffusion with a
stochastic approximate drift is proportional to the expected drift error bound $\tpi(\epsilon)$. 
The bound of \cref{thm:stochastic-approx-grad-distance} thus takes into account how 
the drift error varies with the location of the drift.
Our results match the asymptotic behavior for stochastic gradient Langevin dynamics 
documented in \citet{Teh:2016}: in the limit of the step size 
going to zero, they show that the stochastic gradient has no effect on the equilibrium distribution.

\begin{example}
Suppose $\tY_{t}$ is an Ornstein--Uhlenbeck process with $\ell = d$, the dimensionality of $\tX_{t}$.
That is, for some $\alpha, v > 0$, $\dee\tY_{t} = -\alpha \tY_{t}\dee t + \sqrt{2v}\dee\tW_{t}^{Y}$.
Then the equilibrium distribution of $\tY_{t}$ is that of a Gaussian with 
covariance $\sigma^2I$, where $\sigma^{2} \defined v/\alpha$.
Let $\tb(x, y) = b(x) + y$, so $\EE[\tb(\tX, \tY) \given \tX = x] = b(x)$ and 
hence $d_{\mcW}(\pi, \tpi) = 0$. 
\end{example}

While exponential contractivity is natural and applies in many settings, it is useful 
to have bounds on the Wasserstein distance of approximations when the diffusion 
process mixes more slowly. 
We can prove the analogous guarantee of \cref{thm:approx-grad-distance} when a 
weaker, polynomial contractivity condition is satisfied.  

\begin{assumption}[\textbf{Polynomial contractivity}] \label{asm:polynomial-contractivity}
There exist constants  $C > 0$, $\alpha > 1$, and $\beta > 0$ such that
for all $x, x' \in \reals^d$,
\[
d_{\mcW}(\mu_{x,t}, \mu_{x',t}) \le C\|x - x'\|_{2}(t + \beta)^{-\alpha}.
\]
\end{assumption}

The parameters $\alpha$ and $\beta$ determines how quickly the diffusion converges to equilibrium.
Polynomial contractivity can be certified using, for example, the techniques 
from \citet{Butkovsky:2014} (see also the references therein). 

\bnthm[Error induced by approximate drift, polynomial contractivity] \label{thm:polynomial-contractivity}
Let $\pi$ and $\tpi$ denote the invariant distributions of the 
diffusions in \cref{eqn:diffusion} and \cref{eqn:diffusion-approx}, respectively. 
If the diffusion \cref{eqn:diffusion} is polynomially ergodic with parameters 
$C$, $\alpha$, and $\beta$, the regularity conditions of \cref{asm:regularity} hold, 
and $\|b - \tb\|_{2} \le \eps$, then
\[
d_{\mcW}(\pi, \tpi) \le \frac{C\eps}{(\alpha - 1)\beta^{\alpha - 1}}. \label{eq:poly-approx-bound}
\]
\enthm

\begin{remark}[\emph{Coherency of the error bound}]
The error bound of \cref{eq:poly-approx-bound} has a coherent dependence 
on its parameters, just like  \cref{eq:det-gradient-approx-bound}.
If we change the time scale of the diffusion from $t$ to $s = at$ for some $a > 0$,
the polynomial contractivity constants $C, \alpha$, and $\beta$ become, respectively, $C/a^{\alpha},
\alpha$, and $\beta/a$. 
Making these substitutions and replacing $\epsilon$ by $\epsilon a$,
one can check that the scaling $a$ cancels out in the error bound, so the error is independent 
of how we set the time scale.   
\end{remark} 

\section{Overview of analysis techniques} \label{sec:overview}
We use Stein's method \citep{Stein:1972,Barbour:1990,Ross:2011} to bound the 
Wasserstein distance between $\pi$ and $\tpi$ as a function of a bound 
on $\|b - \tb\|_2$ and the mixing time of $\pi$. 
We describe the analysis ideas for the setting when 
$\|b - \tb\|_2 < \epsilon$ (\cref{thm:approx-grad-distance}); 
the analysis with stochastic drift (\cref{thm:stochastic-approx-grad-distance}) 
or assuming polynomial contractivity (\cref{thm:polynomial-contractivity}) is similar. 
All of the details are in the \suppmat.   

For a diffusion $(X_{t})_{t\ge0}$ with drift $b$,
the corresponding infinitesimal generator satisfies
\[
\gen{b}\phi(x) &= b(x) \cdot \grad \phi(x) + \Delta\phi(x) \label{eq:generator}
\]
for any function $\phi$ that is twice continuously differentiable and vanishing at infinity. 
See, e.g., \citet{Ethier:2009} for an introduction to infinitesimal generators.  
Under quite general conditions, the invariant measure $\pi$ and the 
generator $\gen{b}$ satisfy
\[
\pi(\gen{b}\phi) = 0.
\]
For any measure $\nu$ on $\reals^{d}$ and set of test functions 
$\mcF \subseteq C^{2}(\reals^{d}, \reals)$, we can define the \emph{Stein discrepancy} as: 
\[
\mcS(\nu, \gen{b}, \mcF) 
&\defined \sup_{\phi \in \mcF} |\pi(\gen{b}\phi) - \nu(\gen{b}\phi)|
= \sup_{\phi \in \mcF} |\nu(\gen{b}\phi)|.
\] 
The Stein discrepancy quantifies the difference between $\nu$ and $\pi$ in terms of 
the maximum difference in the expected value of a function (belonging to the transformed 
test class $\{ \gen{b}\phi \given \phi \in \mcF\}$) under these two distributions. 
We can analyze the Stein discrepancy between $\pi$ and $\tpi$ as follows. 
Consider a test set $\mcF$ such that $\|\grad \phi\|_{2} \le 1$ for all 
$\phi \in \mcF$, which is equivalent to having $\|\phi\|_{L} \le 1$. 
We have that 
\[
\mcS(\tpi, \gen{b}, \mcF) 
&= \sup_{\phi \in \mcF} |\tpi(\gen{b}\phi)| 
= \sup_{\phi \in \mcF} |\tpi(\gen{b}\phi - \gen{\tb}\phi)| \\
&= \sup_{\phi \in \mcF} |\tpi(\grad \phi \cdot b - \grad \phi \cdot \tb)|  \\
&\le \sup_{\phi \in \mcF} |\tpi(\|\grad \phi\|_{2} \|b - \tb\|_{2})|
\le \eps,
\]
where we have used the definition of Stein discrepancy, that $\tpi(\gen{\tb}\phi) = 0$, 
the definition of the generator, the Cauchy-Schwartz inequality, that $\|\grad \phi\|_{2} \le 1$,
and the assumption $\|b - \tb\|_{2} \le \eps$. 
It remains to show that the Wasserstein distance satisfies $d_{\mcW}(\pi, \tpi) \le C_{\pi} \mcS(\tpi, \gen{b}, \mcF)$ for 
some constant $C_{\pi}$ that may depend on $\pi$. This would then allow us to conclude that $d_{\mcW}(\pi, \tpi) \le C_{\pi}\epsilon$. 
To obtain $C_{\pi}$, for each 1-Lipschitz function $h$, we construct the solution $u_{h}$ 
to the differential equation
\[
h - \pi(h) = \gen{g}u \label{eq:diff-eq}
\]
and show that $\|\grad u_{h}\|_{2} \le C_{\pi}\|\grad h\|_{2}$.
\opt{XXXXX}{\textcolor{white}{\ref*{eq:diff-eq}}} %
\section{Application: computational--statistical trade-offs} \label{sec:trade-off}

As an application of our results we analyze the behavior of the \emph{unadjusted Langevin Monte Carlo algorithm} 
(ULA)~\citep{Roberts:1996} when approximate gradients of the log-likelihood are used. 
ULA  uses a discretization of the continuous-time Langevin diffusion to approximately sample from 
the invariant distribution of the diffusion.
We prove conditions under which we can obtain more accurate samples by using an approximate drift derived from a Taylor expansion of the exact drift.
 
For the diffusion $(X_{t})_{t \ge 0}$ driven by drift $b$ as defined in \cref{eqn:diffusion} 
and a non-increasing sequence of step sizes $(\gamma_{i})_{i \ge 1}$, the associated ULA Markov chain is 
\opt{arxiv}{
\[
X'_{i+1} &= X'_{i} + \gamma_{i+1}\,b(X'_{i}) + \sqrt{2 \gamma_{i+1}} \xi_{i+1}, & \xi_{i+1} \distiid \distNorm(0, 1).
\label{eqn:ula}
\]
}
\opt{XXXXX}{
\[
X'_{i+1} &= X'_{i} + \gamma_{i+1}\,b(X'_{i}) + \sqrt{2 \gamma_{i+1}} \xi_{i+1}, 
\label{eqn:ula}
\]
where $\xi_{i+1} \distiid \distNorm(0, 1)$.}
Recently, substantial progress has been made in understanding the approximation accuracy 
of ULA~\citep{Dalalyan:2017,Bubeck:2015,Durmus:2016}. 
These analyses show, as a function of the discretization step size $\gamma_i$, 
how quickly the distribution of $X'_i$ converges to the desired target distribution. 

In many big data settings, however, computing $b(X'_i)$ exactly at every step is computationally expensive. 
Given a fixed computational budget, one option is to compute $b(X'_i)$ precisely and run the 
discretized diffusion for a small number of steps to generate samples. 
Alternatively, we could replace $b(X'_i)$ with an approximate drift $\tb(X'_i)$ which is cheaper 
to compute and run the discretized approximate diffusion for a larger number of steps to generate samples. 
While approximating the drift can introduce error, running for more steps can compensate
by sampling from a better mixed chain.
Thus, our objective is to compare the ULA chain using an exact drift initialized at some point 
$x^{\star} \in \reals^d$ to a ULA chain using an approximate drift 
initialized at the same point. 
We denote the exact and approximate drift chains by $X'_{x^{\star},i}$ and $\tX'_{x^{\star},i}$, 
respectively, and denote laws of these chains by $\mu^{\star}_{i}$ and $\tmu^{\star}_{i}$.

For concreteness, we analyze generalized linear models with unnormalized log-densities of the form
\[
\mcL(x) \defined \log \pi_{0}(x) + \sum_{i=1}^{N} \phi_{i}(x \cdot y_{i}),
\]
where $y_{1},\dots, y_{N} \in \reals^{d}$ is the data and $x$ is the parameter. 
In this setting the drift is $b(x) = \grad \mcL(x)$.  
We take $x^{\star} = \argmax_{x} \mcL(x)$ and approximate the drift with a  Taylor expansion around $x^{\star}$:
\opt{arxiv}{
\[
\tb(x) \defined  (H\log\pi_{0})(x^{\star})(x - x^{\star}) + \sum_{i=1}^{N} \phi_{i}''(x^{\star} \cdot y_{i}) y_{i}y_{i}^{\top}(x - x^{\star}),
\label{eqn:approx-ula}
\]
}
\opt{XXXXX}{
\[
\begin{split}
\tb(x) &\defined  (H\log\pi_{0})(x^{\star})(x - x^{\star}) \\
&\phantom{\defined~} + \sum_{i=1}^{N} \phi_{i}''(x^{\star} \cdot y_{i}) y_{i}y_{i}^{\top}(x - x^{\star}),
\label{eqn:approx-ula}
\end{split}
\]
}
where $H$ is the Hessian operator. 
The quadratic approximation of \cref{eqn:approx-ula} basically corresponds to taking a Laplace approximation of the log-likelihood. 
In practice, higher-order Taylor truncation or other approximations can be used, and our analysis can
be extended to quantify the trade-offs in those cases as well. 
Here we focus on the second-order approximation as a simple illustration of the computational--statistical trade-off.   

In order for the Taylor approximation to be well-behaved, we require the prior $\pi_{0}$ and 
link functions $\phi_{i}$ to satisfying some regularity conditions, which are usually easy 
to check in practice.

\begin{assumption}[Concavity, smoothness, and asymptotic behavior of data] \label{asm:phi-simple}
~
\benum
\item The function $\log \pi_{0} \in C^{3}(\reals^{d}, \reals)$ is strongly concave, 
$\|\grad \log \pi_{0}\|_{L} < \infty$,
and $\|H[\partial_{j}\log\pi_{0}]\|_{2} < \infty$ for $j=1,\dots,d$, where $\|\cdot\|_{2}$ denotes the
matrix spectral norm. 
\item For $i=1,\dots,N$, the function $\phi_{i} \in C^{3}(\reals, \reals)$ is strongly concave, 
$\|\phi_{i}'\|_{L} < \infty$, and $ \|\phi_{i}'''\|_{\infty} < \infty$.
\item The data satisfies $\|\sum_{i=1}^{N} y_{i}y_{i}^{\top}\|_{2} = \Theta(N)$.
\eenum
\end{assumption}

We measure computational cost by the number of $d$-dimensional inner products performed. 
Running ULA with the original drift $b$ for $T$ steps costs $TN$ because each step needs 
to compute $x\cdot y_i$ for each of the $N$ $y_i$'s. 
Running ULA with the Taylor approximation $\tb$, we need to compute 
$\sum_{i=1}^{N} \phi_{i}''(x^{\star} \cdot y_{i}) y_{i}y_{i}^{\top}$ once up front, 
which costs $Nd$, and then for each step we just multiply this $d$-by-$d$ matrix with $x - x^{\star}$, 
which costs $d$. 
So the total cost of running approximate ULA for $\tT$ steps is $(\tT + N)d$.

\bnthm[Computational--statistical trade-off for ULA] \label{thm:ula-accuracy-main}
Set the step size $\gamma_i = \gamma_1 i^{-\alpha}$ for fixed $\alpha \in (0,1)$ and 
suppose the ULA of \cref{eqn:ula} is run for $T > d$ steps.
If \cref{asm:phi-simple} holds and $\tT$ is chosen such that the computational cost
of the second-order approximate ULA using drift \cref{eqn:approx-ula} equals that 
of the exact ULA, then $\gamma_1$ may be chosen such that 
\opt{arxiv}{
\[
d_{\mcW}^{2}(\mu^{\star}_{T}, \pi) = \tO\left(\frac{d}{TN}\right)
\qquad \text{and} \qquad
d_{\mcW}^{2}(\tmu^{\star}_{\tT}, \pi) = \tO\left(\frac{d^{2}}{N^{2}T} + \frac{d^{3}}{N^{2}}\right).
\]}
\opt{XXXXX}{
\[
d_{\mcW}^{2}(\mu^{\star}_{T}, \pi) &= \tO\left(\frac{d}{TN}\right) 
\]
and
\[
d_{\mcW}^{2}(\tmu^{\star}_{\tT}, \pi) &= \tO\left(\frac{d^{2}}{N^{2}T} + \frac{d^{3}}{N^{2}}\right).
\]}
\enthm

The ULA procedure of \cref{eqn:ula} has Wasserstein error decreasing like $1/N$
for data size $N$. 
Because approximate ULA can be run for more steps at the same computational cost, 
its error decreases as $1/N^2$. 
Thus, for large $N$ and fixed $T$ and $d$, approximate ULA with drift $\tb$ achieves more 
accurate sampling than ULA with $b$. A conceptual benefit of our results is that we can cleanly 
decompose the final error into the discretization error and the equilibrium bias due to approximate drift. 
Our theorems in \cref{sec:main} quantifies the equilibrium bias, and we can apply existing techniques to bound the discretization error.

\begin{figure*}[tb]
\begin{center}
\begin{subfigure}[b]{\figwidth\textwidth}
    \includegraphics[width=\textwidth]{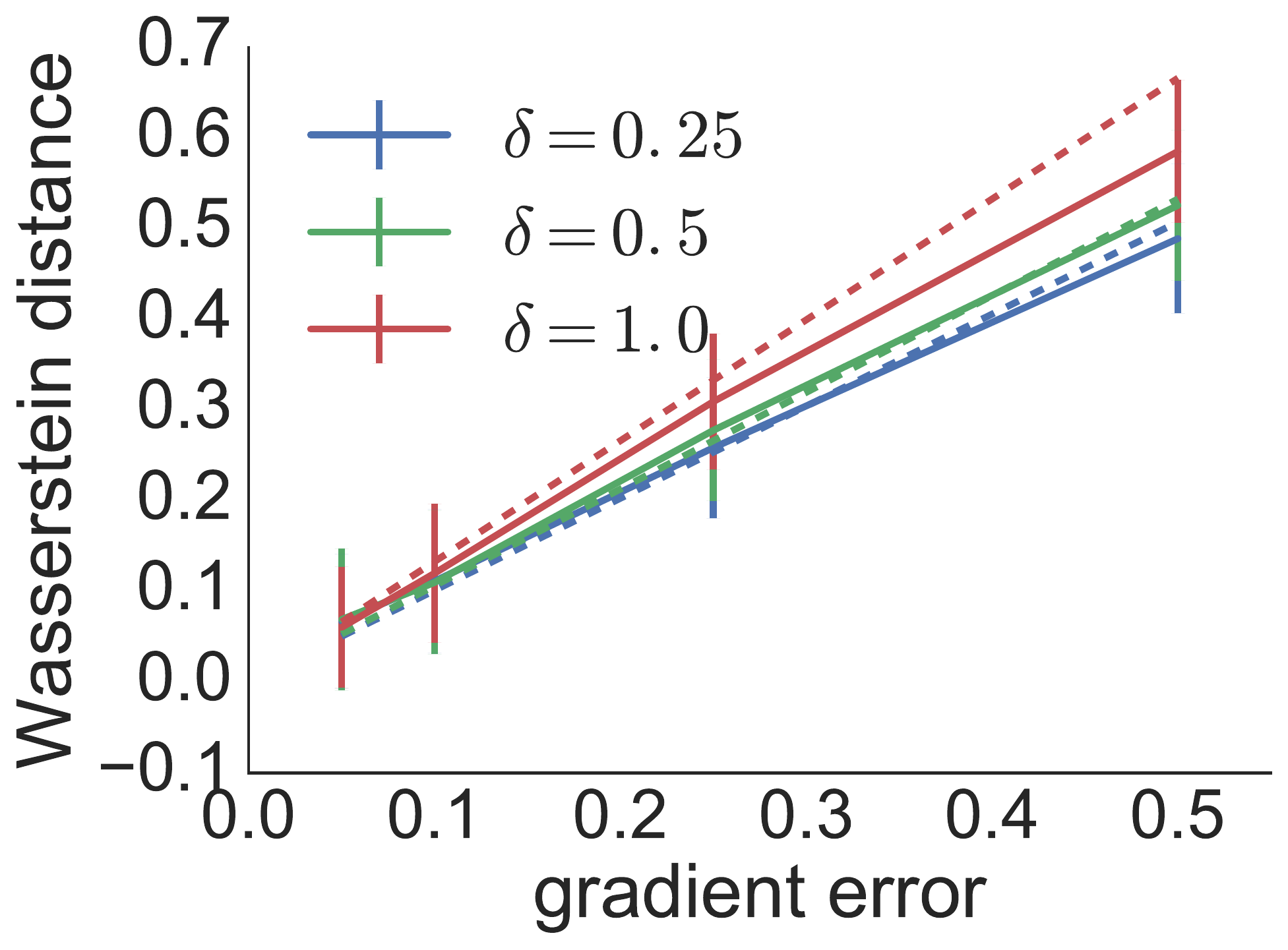}
    \caption{}
    \label{fig:approx-diffusion-small-delta}
\end{subfigure}
\qquad
\begin{subfigure}[b]{\figwidth\textwidth}
    \includegraphics[width=\textwidth]{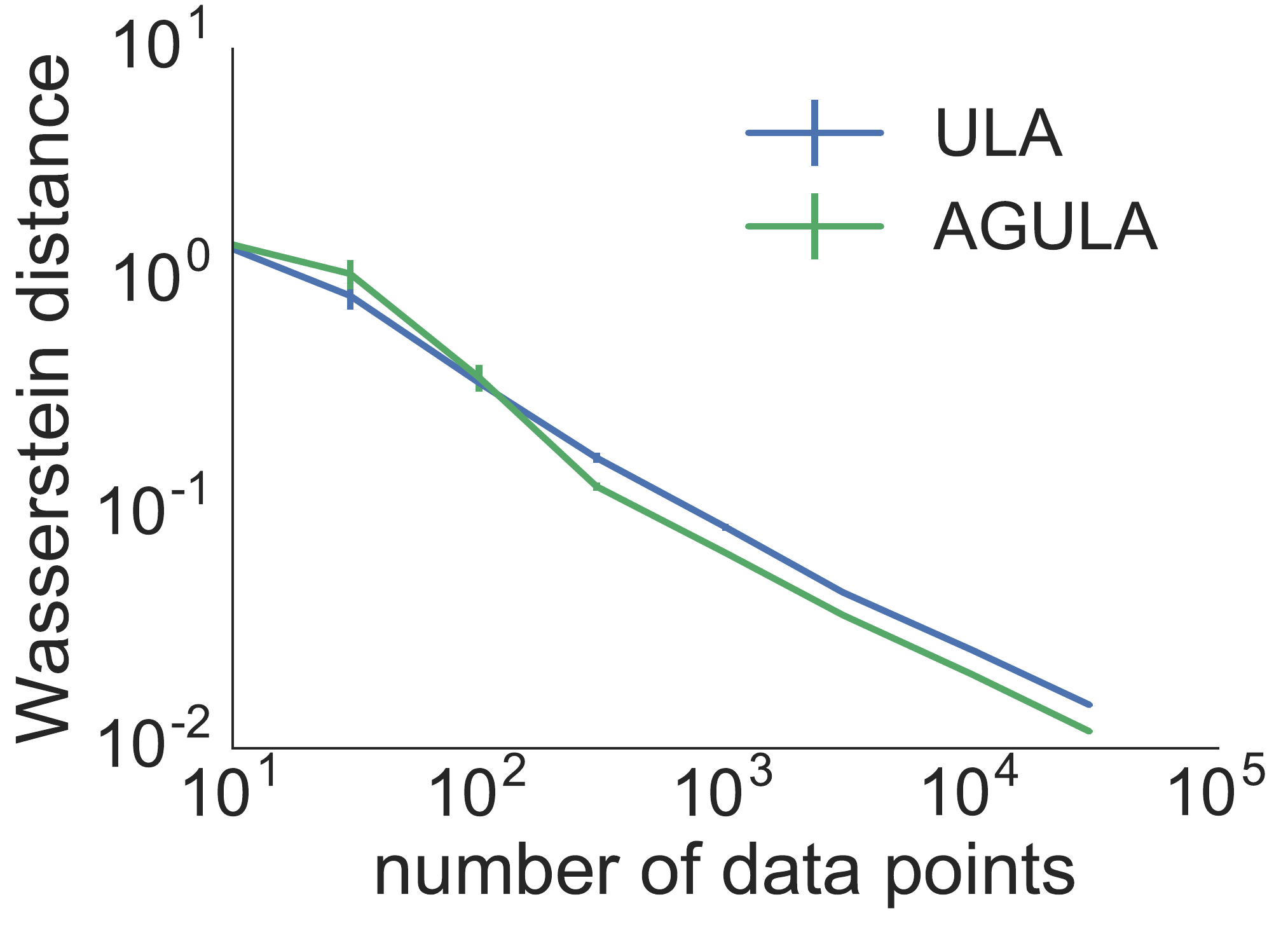}
    \caption{}
    \label{fig:ula}
\end{subfigure}

\opt{XXXXX}{\vspace{-0.2cm}}
\caption{\textbf{(a)} Gradient error $\eps$ versus the Wasserstein distance between $\pi_{\delta}$ 
and $\tpi_{\delta,\eps}$, the stationary distribution of the diffusion with approximate drift 
$\tb_{\delta,\eps}(x) = \grad \log \pi_{\delta}(x) + \eps$.
The solid lines are the simulation results and the dotted lines are the theoretical upper bounds obtained from \cref{thm:approx-grad-distance}.
The simulation results closely match the theoretical bounds and show linear growth in $\eps$, as
predicted by the theory. 
Due to Monte Carlo error the simulation estimates sometimes slightly exceed the theoretical bounds. 
\textbf{(b)} The $y$-axis measures the Wasserstein distance between the true posterior distribution and the finite-time distribution of the exact gradient ULA (ULA) and the approximate gradient ULA (AGULA).
Except for when the number of data points $N < 100$, AGULA shows superior performance, in agreement with the analysis of \cref{thm:ula-accuracy-main}. 
For all experiments the Wasserstein distance was estimated 10 times, each time using 1,000 samples from each distribution.
}
\opt{XXXXX}{\vspace{-0.4cm}}
\label{fig:approx-diffusion}
\end{center}
\end{figure*}

\section{Extension: piecewise deterministic Markov processes}
\label{sec:pdmps}
\vspace{-0.2cm}

We next demonstrate the generality of our techniques by providing a perturbation analysis of
piecewise deterministic Markov processes (PDMPs), which are continuous-time processes that 
are deterministic except at random jump times. 
Originating with the work of \citet{Davis:1984}, there is now a rich literature on 
the ergodic and convergence 
properties of PDMPs~\citep{Costa:2008,Benaim:2012,Fontbona:2012,Azais:2014,Monmarche:2015}.
They have been used to model a range of phenomena including communication networks,
neuronal activity, and biologic population models (see \citep{Azais:2014} and references therein). 
Recently, PDMPs have also been used to design novel MCMC inference schemes.
zig-zag processes (ZZPs)~\citep{Bierkens:2016a,Bierkens:2016b,Bierkens:2016c} are a class of PDMPs that are particularly promising for inference.  ZZPs can be simulated exactly (making Metropolis-Hastings corrections unnecessary) and are non-reversible, which can potentially lead to more efficient sampling~\citep{Mira:2000,Neal:2004}.

Our techniques can be readily applied to analyze the accuracy of approximate PDMPs. 
For concreteness we demonstrate the results for ZZPs in detail and defer the general 
treatment of PDMPs, which includes an idealized version of Hamiltonian Monte Carlo, to the \suppmat.  
The ZZP is defined on the
space $E = \reals^{d} \times \mcB$, where $\mcB \defined \{-1, +1\}^{d}$.
Densities on $\mcB$ are with respect to the counting measure. 

Informally, the behavior of a ZZP can be described as follows. 
The trajectory is $X_t$ and its velocity is $\Theta_t$, so $\frac{d}{dt}X_t = \Theta_t$. 
At random times, a single coordinate of $\Theta_t$ flips signs. 
In between these flips, the velocity is a constant and the trajectory is a straight 
line (hence the name ``zig-zag''). 
The rate at which $\Theta_t$ flips a coordinate is time inhomogeneous. 
The $i$-th component of $\Theta$ switches at rate $\lambda_i(X_t, \Theta_t)$. 
By choosing the switching rates appropriately, the ZZP can be 
made to sample from the desired distribution. 
More precisely, the ZZP $(X_{t}, \Theta_{t})_{t\ge0}$ is determined by the switching rate
$\lambda \in C^{0}(E, \reals_{+}^{d})$ and has  generator
\[
\gen{\lambda}\phi(x, \theta) = \theta \cdot \grad_{x} \phi(x, \theta) + \lambda(x, \theta) \cdot \grad_{\theta} \phi(x, \theta) \label{eq:zzp-generator}
\]
for any sufficiently regular $\phi : E \to \reals$. 
Here $\grad_{x}\phi$ denotes the gradient of $\phi$ with respect to $x$ and  $\grad_{\theta}\phi $ is the discrete differential operator.\footnote{$\grad_{\theta}\phi \defined (\partial_{\theta,1}\phi, \dots, \partial_{\theta,d}\phi)$,
where $\partial_{\theta,i}\phi(x, \theta) \defined \phi(x, R_{i}\theta) - \phi(x, \theta)$
and for $i \in [d]$, the reversal function $R_{i} : \mcB \to \mcB$ is given by 
$
(R_{i}\theta)_{j} \defined 
\begin{cases}
-\theta_{j} & j = i \\
\theta_{j} & j \ne i.
\end{cases}
$} 
Let $(a)^{+} \defined \max(0, a)$ denote the positive part of $a \in \reals$ and $\partial_{i}\phi \defined \D{\phi}{x_{i}}$.
The following result shows how to construct a ZZP with
invariant distribution $\pi$. 

\bnthm[{\citet[Theorem 2.2, Proposition 2.3]{Bierkens:2016a}}]
Suppose $\log \pi \in C^{1}(\reals^{d})$ and $\gamma \in C^{0}(E, \reals^{d}_{+})$
satisfies $\gamma_{i}(x, \theta) = \gamma_{i}(x, R_{i}\theta)$. 
Let 
\[
\lambda_{i}(x, \theta) = (-\theta_{i}\partial_{i}\log \pi(x))^{+} + \gamma_{i}(x, \theta). \label{eq:zzp-lambda}
\]
Then the Markov process with generator $\gen{\lambda}$ has invariant distribution 
$\pi_{E}(\dee x, \theta) = 2^{-d}\pi(\dee x)$. 
\enthm

Analogously to the approximate diffusion setting, we compare 
$(X_{t}, \Theta_{t})_{t\ge0}$ to an approximating ZZP $(\tX_{t}, \tTheta_{t})_{t\ge0}$ 
with switching rate $\tlambda \in C^{0}(E, \reals_{+}^{d})$. 
For example, if $\tpi$ is an approximating density, the approximate 
switching rate could be chosen as
\[
\tlambda_{i}(x, \theta) = (-\theta_{i}\partial_{i}\log \tpi(x))^{+} + \gamma_{i}(x, \theta). \label{eq:zzp-approx-lambda}
\]
To relate the errors in the switching rates to the Wasserstein distance in the final 
distributions, we use the same strategy as before: apply Stein's method to the ZZP generator in \cref{eq:zzp-generator}. 
We rely on ergodicity and regularity conditions that are analogous to those for diffusions. 
We write $(X_{x,\theta,t}, \Theta_{x,\theta,t})$ 
to denote the version of the ZZP satisfying $(X_{x,\theta,0}, \Theta_{x,\theta,0}) = (x, \theta)$
and denote its law by $\mu_{x,\theta,t}$.

\begin{assumption}[\textbf{ZZP polynomial ergodicity}] \label{asm:zzp-polynomial-ergodicity}
There exist constants $C > 0$, $\alpha > 1$, and $\beta > 0$ %
such that for all $x \in \reals^{d}$, $\theta \in \mcB$, and $i \in [d]$,
\[
d_{\mcW}(\mu_{x,\theta,t}, \mu_{x,R_{i}\theta,t}) &\le C(t + \beta)^{-\alpha}.
\]
\end{assumption}
\vspace{-0.2cm}

The ZZP polynomial ergodicity condition is looser than that used for diffusions. 
Indeed, we only need a quantitative bound on the ergodicity constant when the chains 
are started with the same $x$ value.
Together with the fact that $\mcB$ is compact, this
simplifies verification of the condition, which can be done using well-developed coupling
techniques from the PDMP literature~\citep{Benaim:2012,Fontbona:2012,Azais:2014,Monmarche:2015}
as well as more general Lyapunov function-based approaches~\citep{Hairer:2009}.%

Our main result of this section bounds the error in the invariant distributions due 
to errors in the ZZP switching rates. 
It is more natural to measure the error between
$\lambda$ and $\tlambda$ in terms of the $\ell^{1}$ norm.

\bnthm[ZZP error induced by approximate switching rate] \label{thm:zzp-approx-grad-distance}
Assume the ZZP with switching rate $\lambda$ (respectively $\tlambda$) has
invariant distribution $\pi$ (resp.~$\tpi$).
Also assume that $\int_{E} x^{2} \pi(\dee x, \dee \theta) < \infty$ and if
a function $\phi \in C(E, \reals)$ is $\pi$-integrable then it is $\tpi$-integrable. 
If the ZZP with switching rate $\lambda$ is polynomially ergodic with 
constants $C$, $\alpha$, and $\beta$ and $\|\lambda - \tlambda\|_{1} \le \eps$, then
\[
d_{\mcW}(\pi, \tpi) \le \frac{C\eps}{(\alpha - 1)\beta^{\alpha - 1}}. \label{eq:zzp-gradient-approx-bound}
\]
\enthm

\begin{remark}
If the approximate switching rate takes the form of \cref{eq:zzp-approx-lambda}, then
$\|\grad \log \pi - \grad \log \tpi\|_{1} \le \eps$ implies $\|\lambda - \tlambda\|_{1} \le \eps$.
\end{remark}

\opt{XXXXX}{\vspace{-0.2cm}}
\section{Experiments}
\label{sec:experiments}
\opt{XXXXX}{\vspace{-0.2cm}}
We used numerical experiments to investigate whether our bounds capture the true behavior of approximate diffusions and
their discretizations. 

\textbf{Approximate Diffusions.}
For our theoretical results to be a useful guide in practice, we would like the Wasserstein bounds 
to be reasonably tight and have the correct scaling in the problem parameters (e.g., in $\|b - \tb\|_{2}$). 
To test our main result concerning the error induced from using an approximate drift (\cref{thm:approx-grad-distance}),
we consider mixtures of two Gaussian densities of the form
\[
\pi_{\delta}(x) = \frac{1}{2(2\pi)^{d/2}}\left(e^{-\|x-\delta/2\|_{2}^{2}/2} + e^{-\|x+\delta/2\|_{2}^{2}/2}\right),
\]
where $\delta \in \reals^{d}$ parameterizes the difference between the means of the Gaussians. 
If $\|\delta\|_{2} < 2$, then $\pi_{\delta}$ is $(1-\|\delta\|_{2}/4)$-strongly log-concave;
if $\|\delta\|_{2} = 2$, then $\pi_{\delta}$ is log-concave; and if $\|\delta\|_{2} > 2$, then
$\pi_{\delta}$ is not log-concave, but is log-concave in the tails. 
Thus, for all choices of $\delta$, the diffusion with drift $b_{\delta}(x) \defined \grad \log \pi_{\delta}(x)$ is 
exponentially ergodic.
Importantly, this class of Gaussian mixtures allows us to investigate a range of practical regimes, 
from strongly unimodal to highly multi-modal distributions. 
For $d=1$ and a variety of choices of $\delta$, we generated 1,000 samples from the target distribution
$\pi_{\delta}$ (which is the stationary distribution of a diffusion with drift 
$b_{\delta}(x)$) and from $\tpi_{\delta,\eps}$ (which is the stationary distribution of the approximate
diffusion with drift $\tb_{\delta,\eps}(x) \defined b_{\delta}(x) + \eps$) for $\eps = 0.05, 0.1, 0.25, 0.5$. 
We then calculated the Wasserstein distance between the empirical distribution of the target and
the empirical distribution of each approximation.
\cref{fig:approx-diffusion-small-delta} shows the empirical Wasserstein distance (solid lines) 
for $\delta = 0.25, 0.5, 1.0$ along with the corresponding theoretical bounds from 
\cref{thm:approx-grad-distance} (dotted lines). 
The two are in close agreement. We also investigated larger distances for  
$\delta = 1.0, 2.0, 3.0$. Here the exponential contractivity constants that can be derived from \citet{Eberle:2015} are rather loose. 
Importantly, however, for all values of $\delta$ considered, the Wasserstein distance grows 
linearly in $\epsilon$, as predicted by our theory.
Results for $d > 1$ show similar linear behavior in $\epsilon$, though we omit the plots. %

\textbf{Computational--statistical trade-off.}
We illustrate the computational--statistical trade-off of \cref{thm:ula-accuracy-main}
in the case of logistic regression. This corresponds to $\phi_{i}(t) = \phi_{lr}(t) \defined -\log(1 + e^{-t})$. 
We generate data $y_{1},y_{2},\dots$ according to the following process:
\[
z_{i} &\dist \distBern(.5),  &
\zeta_{i} &\dist \distNorm(\mu_{z_{i}}, I), &
y_{i} &= (2z_{i} - 1)\zeta_{i},
\]
where $\mu_{0} = (0, 0, 1, 1)$ and $\mu_{1} = (1, 1, 0, 0)$. 
We restrict the domain $\mcX$ to a ball of radius 3, 
$\mcX = \{ x \in \reals^{4} \given \|x\|_{2} \le 3\}$, and
add a projection step to the ULA algorithm~\citep{Bubeck:2015}, 
replacing $Z'_{i}$ with $\argmin_{z \in \mcX} \|Z'_{i} - z\|_{2}$.
While \cref{thm:ula-accuracy-main} assumes $\mcX = \reals^4$, 
the numerical results here on the bounded domain still illustrate our key point: for the same computational budget, computing fast approximate gradients and 
running the ULA chain for longer can produce a better sampler. 
\cref{fig:ula} shows that except for very small $N$, the approximate gradient ULA (AGULA), which uses the approximation in \cref{eqn:approx-ula},
produces better performance than exact gradient ULA (ULA) with the same budget. 
For each data-set size ($N$), the true posterior distribution was estimated by running an adaptive Metropolis-Hastings (MH)
sampler for 100,000 iterations. 
ULA and AGULA were each run 1,000 times to empirically estimate the
approximate posteriors.
We then calculated the Wasserstein distance between the ULA and AGULA empirical distributions and
the empirical distribution obtained from the MH sampler. 

\opt{XXXXX}{\vspace{-0.2cm}}
\section{Discussion}
\opt{XXXXX}{\vspace{-0.2cm}}
\textbf{Related Work.}
Recent theoretical work on scalable MCMC algorithms has yielded numerous insights into the
regimes in which such methods produce computational gains~\citep{Pillai:2014,Rudolf:2015,Alquier:2016,Johndrow:2015,Johndrow:2016}. 
Many of these works focused on approximate Metropolis-Hastings algorithms,
rather than gradient-based MCMC.
Moreover, the results in these papers are for discrete chains, whereas our results also apply to 
continuous diffusions as well as other continuous-time Markov processes such as the zig-zag process. 
Perhaps the closest to our work is that of \citet{Rudolf:2015} and \citet{Gorham:2016}. The former studies general
perturbations of Markov chains and includes an application to stochastic Langevin dynamics.
They also rely on a Wasserstein contraction condition, like our \cref{asm:exponential-contractivity}, 
in conjunction with a Lyapunov condition on the perturbed chain.
However, our more specialized analysis is particularly transparent and leads to tighter 
bounds in terms of the contraction constant $\rho$: the bound of \citet{Rudolf:2015} is proportional
to $(1-\rho)^{-1}$ whereas our bound is proportional to $-(\log \rho)^{-1}$. 
Another advantage of our approach is that our results are more straightforward to apply since 
we do not need to directly analyze the Lyapunov potential and the perturbation ratios as in \citet{Rudolf:2015}.
Our techniques also apply to the weaker polynomial contraction setting. 
\citet{Gorham:2016} have results of similar flavor to ours and also rely on Stein's method, but their assumptions and
target use cases differ from ours. 
Our results in \cref{sec:trade-off}, which apply when ULA is used with a deterministic approximation to the drift, complement 
the work of \citet{Teh:2016} and \citet{Vollmer:2016}, which provides (non-)asymptotic
analysis when the drift is approximated stochastically at each iteration.  

\textbf{Conclusion.} 
We have established general results on the accuracy of diffusions with approximate drifts. 
As an application, we show how this framework can quantify the computational--statistical 
trade-off in approximate gradient ULA. The example in \cref{sec:experiments} illustrates 
how the log-concavity constant can be estimated in practice and how theory provides reasonably precise error bounds.
We expect our general framework to have many further applications. 
In particular, an interesting direction is to extend our framework 
to analyze the trade-offs in subsampling Hamiltonian Monte Carlo algorithms and stochastic Langevin dynamics.

\subsubsection*{Acknowledgments}
Thanks to Natesh Pillai for helpful discussions and to Trevor Campbell for feedback on an earlier draft. 
Thanks to Ari Pakman for pointing out some typos and to Nick Whiteley for noticing \cref{thm:ula-accuracy-main}
was missing a necessary assumption. 
JHH is supported by the U.S. Government under FA9550-11-C-0028 and awarded by the DoD, Air Force
Office of Scientific Research, National Defense Science and Engineering Graduate
(NDSEG) Fellowship, 32 CFR 168a.

\appendix

\renewcommand*{\theassumption}{\Alph{section}.\arabic{assumption}}

\section{Exponential contractivity}

A natural generalization of the strong concavity case is to assume that $\log \pi$ 
is strongly concave for $x$ and $x'$ far apart and that $\log \pi$ has ``bounded convexity'' 
when $x$ and $x'$ are close together.  
It turns out that in such cases \cref{asm:exponential-contractivity} still holds.
More formally, the following assumption can be used even when the drift is not a gradient. 
For $f : \mcX \to \reals^{d}$ and $r > 0$, let 
\[ 
\kappa(r) \defined \inf\left\{ -2\frac{(f(x) - f(x')) \cdot (x - x')}{r^{2}} : x, x' \in \mcX, \|x - x'\|_{2} = r\right\}.
\] 
Define the constant $R_{0} = \inf\{R \ge 0 : \kappa(r) \ge 0~\forall r \ge R\}$. 

\begin{assumption}[Strongly log-concave tails] \label{asm:strongly-log-concave-in-tails}
For the function $f \in C^{1}(\mcX, \reals^{d})$, there exist constants $R, \ell \in [0,\infty)$
and  $k \in (0, \infty)$  such that 
\[
\begin{aligned}
\kappa(r) \ge -\ell \text{ for all } r \le R \text{ and } \kappa(r) \ge k \text{ for all } r > R.
\end{aligned}  \label{eq:strongly-log-concave-in-tails}
\]
Furthermore, $\kappa(r)$ is continuous and $\int_{0}^{1} r\kappa(r)^{-}\dee r < \infty$.
\end{assumption}

\bnthm[\citet{Eberle:2015,Wang:2016}] \label{thm:exponential-contractivity-general}
If \cref{asm:strongly-log-concave-in-tails} holds for $f = b$ then \cref{asm:exponential-contractivity}
holds for
\[
C &= \exp\left(\frac{1}{4}\int_{0}^{R_{0}} r \kappa(r)\dee s\right) \\
\frac{1}{\log(1/\rho)} &\le 
\begin{cases}
\frac{3e}{2} \max(R^{2}, 8k^{-1}) & \text{if } \ell R_{0}^{2} \le 8 \\
8\sqrt{2\pi}R^{-1}\ell^{-1/2}(\ell^{-1}+k^{-1})e^{\ell R^{2}/8} + 32R^{-2}k^{-2} & \text{otherwise.}
\end{cases}
\]
\enthm

For detailed calculations for the case of a mixture of Gaussians model, see~\citet{Gorham:2016}.

\section{Proofs of the main results in Section~\ref{sec:main}}

We state all our results in the more general case of a diffusion on
a convex space $\mcX \subseteq \reals^{d}$. 
We begin with some additional definitions. 
Any set $\mcG \subseteq C(\mcX)$ defines an \emph{integral probability metric} (IPM)
\[
d_{\mcG}(\mu, \nu) &= \sup_{\phi \in \mcG} |\mu(\phi) - \nu(\phi)|, 
\]
where $\mu$ and $\nu$ are measures on $\mcX$. 
The \emph{Wasserstein metric} $d_{\mcW}$ corresponds to $\mcW \defined \{ \phi \in C(\mcX) \given \|\phi\|_{L} \le 1\}$.
The set $\mcH \defined \{ \phi \in C^{1}(\mcX) \given \|h\|_{L} \le 1 \}$ will be used to define an IPM $d_{\mcH}$. 
For a set $\mcZ \subseteq \reals^{n}$, we use $\partial \mcZ$ to denote the boundary of $\mcZ$.

Suppose $\|b - \tb\|_2 \le \epsilon$.
We first state several standard properties of the Wasserstein metric and invariant measures of diffusions. 
The proofs are included here for completeness. 

\bnlem \label{lem:dH = dW}
For any $\mu, \nu \in \mcP(\mcX)$,
$d_{\mcH}(\mu, \nu) = d_{\mcW}(\mu, \nu)$. 
\enlem

\bprf[Proof sketch]
The result follows since any Lipschitz function is continuous and \aev-differentiable, and
continuously differentiable functions are dense in the class of continuous and 
\aev-differentiable functions. 
\eprf

We use the notation $(X_t)_{t\ge 0} \dist \distDiff(b, \Sigma) $ if $X_t$ is a diffusion defined by 
\[
\dee X_{t} = b(X_t)\,\dee t + \Sigma \dee W_{t} - n_t L(\dee t).
\]
 A diffusion $X_{t}$ is said to be \emph{strong Feller} if its semigroup operator
$(\pi_{t}\phi)(x) \defined \EE[\phi(X_{x,t})]$, $\phi \in C(\mcX)$, satisfies the property
that for all bounded $\phi$, $\pi_{t}\phi$ is bounded and continuous. 

\bnprop \label{prop:pi-diffusion-properties}
Assume \cref{asm:regularity}(1) holds and let $(X_{t})_{t\ge 0} \dist \distDiff(b, I)$.
Then for each $x \in \mcX$, $X_{x,t}$ has the invariant density $\pi$ and is strong Feller.
\enprop
\bprf
The existence of the diffusions follows from \citet[Theorem 4.1]{Tanaka:1979},
the strong Feller property follows from \citet[Ch.~8, Theorems 1.5 \& 1.6]{Ethier:2009},
and the fact that $\pi$ is the unique stationary measure follows since $\gen{b}^{*}\pi = 0$. 
\eprf

By the same proof as \cref{prop:pi-diffusion-properties}, we have

\bnprop[Diffusion properties] \label{prop:diffusion-properties}
For $f \in C^{0}(\mcX, \reals^{d})$ with $\|f\|_{L} < \infty$, 
the diffusion $(X_{t})_{t\ge 0} \dist \distDiff(f, I)$ exists and has an invariant distribution $\pi_{f}$. 
\enprop

\bnprop[Expectation of the generator] \label{prop:expected-b-generator-zero}
For $f \in C^{0}(\mcX, \reals^{d})$, let the diffusion $(X_{t})_{t\ge 0} \dist \distDiff(f, I)$
have invariant density $\pi_{f}$ and assume that linear functions are $\pi_{f}$-integrable.
Then for all $\phi \in C^{2}(\mcX)$ such that $\|\phi\|_{L} < \infty$ and $\gen{f}\phi$ is 
$\pi_{f}$-integrable, $\pi_{f}(\gen{f}\phi) = 0$.
\enprop

\bprf
Let $P_{t}$ be the semigroup operator associated with $(X_{t})_{t \ge 0}$:
\[
(P_t\phi)(x) = \EE[\phi(X_{x,t})].
\]
Since by hypothesis linear functions are $\pi_{f}$-integrable and $\phi$ is Lipschitz, $\phi$ is 
$\pi_{f}$-integrable. 
Thus, $P_{t}\phi$ is $\pi_{f}$-integrable and by the definition of an invariant measure
(see \citep[Definition 1.2.1]{Bakry:2014} and subsequent discussion),
\[
\pi_{f}(P_{t}\phi) = \pi_{f}\phi. \label{eq:invariant-measure-condition}
\]
Using the fact that $\partial_{t}P_{t} = P_{t}\gen{f}$~\citep[Eq.~(1.4.1)]{Bakry:2014}, 
differentiating both size of \cref{eq:invariant-measure-condition}, 
applying dominated convergence, and using the hypothesis that $\gen{f}\phi$ is 
$\pi_{f}$-integrable yields
\[
0 
= \partial_{t}\pi_{f}(P_{t}\phi) 
= \pi_{f}(\partial_{t}P_{t}\phi)
= \pi_{f}(P_{t}\gen{f}\phi)
= \pi_{f}(\gen{f}\phi).
\]
\eprf

We next show that the solution to \cref{eq:diff-eq} is 
Lipschitz continuous with a Lipschitz constant depending on the mixing properties
of the diffusion associated with the generator. 

\bnprop[Differential equation solution properties] \label{prop:integral-transformation}
If Assumptions \ref{asm:exponential-contractivity} and \ref{asm:regularity}(1) hold, then 
for any $h \in C^{1}(\mcX)$ with $\|h\|_{L} < \infty$, the function
\[
u_{h}(x) \defined \int_{0}^{\infty}(\pi(h) - \EE[h(X_{x,t})])\,\dee t
\]
exists and satisfies
\[
\|u_{h}\|_{L} &\le \frac{C}{\log(1/\rho)} \|h\|_{L} \label{eq:u-h-lipschitz}\\
(\gen{b}u_{h})(x) &= h(x) - \pi(h). \label{eq:stein-equation}
\] 
\enprop

\bprf
We follow the approach of \citet{Mackey:2015}. 
By \cref{asm:exponential-contractivity} and the definition of Wasserstein distance, 
we have that there is a coupling between $X_{x,t}$ and $X_{x',t}$ such that
\[
\EE[\|X_{x,t} - X_{x',t}\|_{2}] \le C\|x - x'\|_{2}\rho^{t}.
\]
The function $u_{h}$ is well-defined since for any $x \in \mcX$,
\[
\int_{0}^{\infty}|\pi(h) - \EE[h(X_{x,t})]| \,\dee t 
&= \int_{0}^{\infty}\left|\int_{\mcX}\left(\EE[h(X_{x',t})] - \EE[h(X_{x,t})]\right)\pi(x')\,\dee x'\right|\,\dee t \\
&\le \sup_{z \in \mcX} \|\grad h(z)\|_{2}\int_{0}^{\infty}\int_{\mcX} \EE[\|X_{x,t} - X_{x',t}\|_{2}]\pi(x')\,\dee x'\,\dee t  \\
& = \sup_{z \in \mcX} \|\grad h(z)\|_{2}\int_{0}^{\infty}\int_{\mcX} \|x - x'\|_2 C\rho^t\pi(x')\,\dee x'\,\dee t\\
&\le \|h\|_{L}\,\EE_{X \dist \pi}[\|x - X\|_{2}]\int_{0}^{\infty}C\rho^{t}\,\dee t  \\
&< \infty,
\]
where the first line uses the property that $\pi(h) = \int_{\mcX} \EE[h(X_{x',t})]\pi(x')\dee x'$ and 
the final inequality follows from \cref{asm:regularity}(1) and the assumption that $0 < \rho < 1$. 
Furthermore, $u_{h}$ has bounded Lipschitz constant since for any $x, x' \in \mcX$, 
\[
|u_{h}(x) - u_{h}(x')|
&= \left|\int_{0}^{\infty}\EE[h(X_{x,t}) - h(X_{x',t})]\,\dee t\right| \\
&\le \sup_{z \in \mcX} \|\grad h(z)\|_{2}\int_{0}^{\infty}\EE[\|X_{x,t} - X_{x',t}\|_{2}]\,\dee t \\
&\le \|h\|_{L} \|x - x'\|_{2}  \int_{0}^{\infty} C \rho^{t}\,\dee t \\
&= \frac{C\|h\|_{L}}{\log(1/\rho)} \|x - x'\|_{2}.
\]
Finally, we show that $(\gen{b}u_{h})(x) = h(x) - \pi(h)$.
Recall that for $h \in C(\mcX)$, the semigroup operator is given by
$(\pi_{t}h)(x) = \EE[h(X_{x,t})]$.
Since $X_{x,t}$ is strong Feller for all $x \in \mcX$ by \cref{prop:pi-diffusion-properties},
for all $t \ge 0$, its generator satisfies \citep[Ch.~1, Proposition~1.5]{Ethier:2009}
\[
h - \pi_{t}h = \gen{b}\int_{0}^{t}(\pi(h) - \pi_{s}h)\,\dee s. \label{eq:generator-formula}
\]
Hence,
\[
\begin{split}
\lefteqn{|h(x) - \pi(h) - [h(x) - (\pi_{t}h)(x)]|} \\
&= \left|\int_{\mcX}\EE[h(X_{x,t})] - \EE[h(X_{x',t})] \pi(x')\,\dee x'\right|  \\
&\le \sup_{z \in \mcX} \|\grad h(z)\|_{2}\int_{\mcX} \EE[\|X_{x',t} - X_{x,t}\|_{2}]\pi(x')\,\dee x' \\
&\le \|h\|_{L}\,\EE_{X \dist \pi}[\|x - X\|_{2}]C\rho^{t}.
\end{split}
\]
Thus, conclude that the left-hand side of \cref{eq:generator-formula} converges
pointwise to $h(x) - \pi(h)$ as $t \to \infty$.
Since $\gen{b}$ is closed~\citep[Ch.~1, Proposition~1.6]{Ethier:2009}, the right-hand side of
\cref{eq:generator-formula} limits to $\gen{b} u_{h}$. 
Hence, $u_{h}$ solves \cref{eq:stein-equation}. 
\eprf

We can now prove the main result bounding the Wasserstein distance between the invariant distributions
of the original and perturbed diffusions.

\bprf[Proof of \cref{thm:approx-grad-distance}]
~
By \cref{prop:diffusion-properties} and \cref{asm:regularity}, 
the hypotheses of \cref{prop:expected-b-generator-zero} hold for $f = \tb$. 
Let $\mcF \defined \{ u_{h} \given h \in \mcH \}$. 
Then
\[
d_{\mcW}(\pi, \tpi) 
&= \sup_{h \in \mcH} |\pi(h) - \tpi(h)| \quad \text{by definition and \cref{asm:regularity}} \\
&= \sup_{h \in \mcH} |\pi(\gen{b}u_h) - \tpi(\gen{b}u_{h})| \quad \text{by \cref{eq:stein-equation}} \\
&= \sup_{h \in \mcH} |\tpi(\gen{b}u_{h})|   \quad \text{by \cref{prop:expected-b-generator-zero}} \\
&= \sup_{u \in \mcF} |\tpi(\gen{b}u)|   \quad \text{by definition of $\mcF$} \\
&= \sup_{u \in \mcF} |\tpi(\gen{b}u - \gen{\tb}u)|  \quad \text{by \cref{prop:expected-b-generator-zero}}\\
&= \sup_{u \in \mcF} |\tpi(\grad u \cdot b - \grad u \cdot \tb)| \qquad \text{by definition of $\mcA_{b}$} \\
&\le \sup_{u \in \mcF} |\tpi(\|\grad u\|_{2} \|b - \tb)\|_{2})| \\
&\le \frac{C\eps}{\log(1/\rho)} \qquad \text{by \cref{eq:u-h-lipschitz} and $\|b - \tb\|_{2} \le \eps$}.
\]
\eprf

A similar analysis can be used to bound the Wasserstein distance between $\pi$ and $\tpi$ when the approximate drift $\tb$ is itself stochastic. 

\bprf[Proof of \cref{thm:stochastic-approx-grad-distance}]
We will need to consider the joint diffusions $Z_{t} = (X_{t}, Y_{t})$ and $\tZ_{t} = (\tX_{t}, \tY_{t})$
on $\mcZ \defined \mcX \times \reals^{d}$, where
\[
\dee Z_{t}  &= (b(X_{t}), b_{aux}(Y_{t}))\,\dee t + (\sqrt{2} \dee W^{X}_{t}, \Sigma\,\dee W^{Y}_{t}) - n_{t}L(\dee t) \\
\dee\tZ_{t} &= (\tb(\tX_{t}, \tY_{t}), b_{aux}(\tY_{t}))\,\dee t + (\sqrt{2} \dee \tW^{X}_{t}, \Sigma\,\dee \tW^{Y}_{t}) - n_{t}\tL(\dee t).
\]
Notice that $X_{t}$ and $Y_{t}$ are independent and the invariant distribution of $X_{t}$ is $\pi$. 
Let $\pi_{Z}$ and $\tpi_{Z}$ be the invariant distributions of $Z_{t}$ and $\tZ_{t}$, respectively.
Also note that the generators for $Z_{t}$ and $\tZ_{t}$ are, respectively,
\[
\gen{Z}\phi(z) &= \grad \phi \cdot (b(x), b_{aux}(y)) + \Delta\phi_{x}(z) + \Sigma^{\top}\Sigma : H \phi_{y}(z) \\
\gen{\tZ}\phi(z) &= \grad \phi \cdot (\tb(x, y), b_{aux}(y)) + \Delta\phi_{x}(z) + \Sigma^{\top}\Sigma : H \phi_{y}(z).
\]
where $H$ is the Hessian operator.  
 
By \cref{prop:diffusion-properties} and \ref{asm:regularity}, the hypotheses of \cref{prop:expected-b-generator-zero} hold for $f(x,y) = (\tb(x, y), b_{aux}(y))$.
Let $\mcH_{Z} \defined \{ h \in C^{1}(\mcZ) \given \|h\|_{L} \le 1\}$ 
and $\mcF_{Z} \defined \{ u_{h} \given h \in \mcH_{Z}\}$. 
Also, for $z = (x,y) \in \mcZ$, let $\id_{Y}(z) = y$. 
Then, by reasoning analogous to that in the proof of \cref{thm:approx-grad-distance},
\[
d_{\mcW}(\pi, \tpi) 
&\le d_{\mcW}(\pi_{Z}, \tpi_{Z}) \\
&= \sup_{h \in \mcH_{Z}} |\pi_{Z}(h) - \tpi_{Z}(h)| \\
&= \sup_{u \in \mcF_{Z}} |\tpi_{Z}(\gen{Z}u - \gen{\tZ}u)| \\
&= \sup_{u \in \mcF_{Z}} |\tpi_{Z}(\grad u \cdot b - \grad u \cdot \tb)| \\
&= \sup_{u \in \mcF_{Z}} |\EE[\grad u(\tX, \tY) \cdot \EE[b(\tX) - \tb(\tX, \tY) \given \tX]]| \\
&\le \sup_{u \in \mcF_{Z}} |\EE[\|\grad u(\tX, \tY)\|_{2}\|\EE[b(\tX) - \tb(\tX, \tY) \given \tX]\|_{2}]| \\
&\le \frac{C\,\tpi(\eps)}{\log(1/\rho)}.
\]
\eprf

\bprf[Proof of \cref{thm:polynomial-contractivity}]
The proof is very similar to that of \cref{thm:approx-grad-distance}, the only difference is in the Lipshitz coefficient of the differential equation solution $u_h(x)$ in \ref{prop:integral-transformation}. Using polynomial contractivity, we have
\[
|u_{h}(x) - u_{h}(x')|
&= \left|\int_{0}^{\infty}\EE[h(X_{x,t}) - h(X_{x',t})]\,\dee t\right| \\
&\le \sup_{z \in \mcX} \|\grad h(z)\|_{2}\int_{0}^{\infty}\EE[\|X_{x,t} - X_{x',t}\|_{2}]\,\dee t \\
&\le \|h\|_{L} \|x - x'\|_{2}  \int_{0}^{\infty} C (t+\beta)^{-\alpha}\,\dee t \\
&= \frac{C\|h\|_{L}}{(\alpha - 1)\beta^{\alpha-1}} \|x - x'\|_{2}.
\]
Plugging in this Lipschitz constant, we have 
\[
d_{\mcW}(\pi, \tpi) \leq \frac{C\epsilon}{(\alpha-1)\beta^{\alpha-1}}.
\]
\eprf 

\section{Checking the Integrability Condition}

The following result gives checkable conditions under which \cref{asm:regularity}(3) holds.
Let $\ball_{R} \defined \{ x \in \reals^{d} \given \|x\|_{2} \le R\}$.

\bnprop[Ensuring integrability] \label{prop:checking-integrability}
\cref{asm:regularity}(3) is satisfied if $b = \grad \log \pi$, $\tb = \grad \log \tpi$, 
$\|b - \tb\|_{2} \le \eps$, and either
\benum
\item there exist constants $R > 0, B > 0, \delta > 0$ such that for all
$x \in \mcX \setminus \ball_{R}$, $\|b(x) - \tb(x)\|_{2} \le B/\|x\|_{2}^{1+\delta}$; or
\item there exists a constant $R > 0$ such that for all $x \in \mcX \setminus \ball_{R}$
$x \cdot (b(x) - \tb(x)) \ge 0$.
\eenum
\enprop 
\bprf
For case (1), first we note that since $\int_{\mcX}(\pi(x) - \tpi(x))\,\dee x = 0$,
by the (generalized) intermediate value theorem, there exists $x^{*} \in \mcX$
such that $\pi(x^{*}) - \tpi(x^{*}) = 0$, and hence 
$\log \pi(x^{*}) - \log \tpi(x^{*}) = 0$.
Let $p[x^{*},x]$ be any path from $x^{*}$ to $x$. 
By the fundamental theorem of calculus for line integrals,
\[
|\log \pi(x) - \log \tpi(x)|
&= \Big|\log \tpi(x^{*}) - \log \pi(x^{*}) + \int_{\gamma[x^{*},x]}(b(r) - \tb(r))\cdot \dee r\Big| \\
&= \Big|\int_{\gamma[x^{*},x]}(b(r) - \tb(r))\cdot r'(t) \,\dee t\Big| \\
&\le \int_{\gamma[x^{*},x]}\|b(r) - \tb(r)\|_{2}\|r'(t)\|_{2} \,\dee t.
\]
First consider $x \in \mcX \cap \ball_{R}$. 
Choosing $p[x^{*},x]$ to be the linear path $\gamma[x^{*},x]$, we have
\[
|\log\pi(x) - \log \tpi(x)|  
&\le \eps\int_{\gamma[x^{*},x]} \|r'(t)\|_{2} \,\dee t \\
&= \eps \|x - x^{*}\|_{2} \\
&\le (R + \ell^{*})\eps,  \label{eq:log-density-in-ball-diff}
\]
where $\ell^{*} \defined \|x^{*}\|_{2}$.

Next consider $x \in \mcX \setminus \ball_{R}$. 
Let $\ell \defined \|x\|_{2}$ and $x' = \frac{R}{\ell}x$. 
Choose $p[x^{*},x]$ to consist of the concatenation of the linear paths $\gamma[x^{*},0]$,
$\gamma[0,x']$, and $\gamma[x',0]$, so
\[
\begin{split}
\lefteqn{\int_{p[x^{*}, x]} \|b(r)-\tb(r)\|_{2}\|r'(t)\|_{2} \,\dee t} \\
&= \int_{\gamma[x^{*},0]} \|b(r)-\tb(r)\|_{2}\|r'(t)\|_{2} \,\dee t + \int_{\gamma[0,x']} \|b(r)-\tb(r)\|_{2}\|r'(t)\|_{2} \,\dee t \\
&\phantom{=~} + \int_{\gamma[x',x]} \|b(r)-\tb(r)\|_{2}\|r'(t)\|_{2} \,\dee t.
\end{split}
\]
Now, we bound each term:
\[
\int_{\gamma[x^{*},0]} \|b(r)-\tb(r)\|_{2}\|r'(t)\|_{2} \,\dee t  &\le \ell^{*}\eps \\
\int_{\gamma[0,x']} \|b(r)-\tb(r)\|_{2}\|r'(t)\|_{2} \,\dee t &\le R\eps \\
\int_{\gamma[x',x]} \|b(r)-\tb(r)\|_{2}\|r'(t)\|_{2} \,\dee t
&\le (\ell - R)B\int_{0}^{1} \frac{1}{(R + (\ell - R)t)^{1+\delta}} \\
&= (\ell - R)B\left[ \frac{1}{(\ell - R)R^{\delta}} -  \frac{1}{(\ell - R)\ell^{\delta}}\right] \\
&\le \frac{B}{R^{\delta}}. 
\]
It follows that there exists a constant $\tB > 0$ such that 
for all $x \in \mcX$, $|\log\pi(x) - \log \tpi(x)| < \tB$. 
Hence $\tB^{-1}\pi < \tpi < \tB\pi$,
so $\phi$ is $\pi$-integrable if and only if it is $\tpi$-integrable. 

Case (2) requires a slightly more delicate argument. 
Let $x^{*}$ and $\ell^{*}$ be the same as in case (1).
For $x \in \mcX \cap \ball_{R}$, it follows from \cref{eq:log-density-in-ball-diff} that
\[
\log\pi(x) - \log \tpi(x) \ge -(R + \ell^{*})\eps. 
\]
For $x \in \mcX \setminus \ball_{R}$, arguing as above yields
\[
\log\pi(x) - \log \tpi(x)
&= \int_{p[x^{*},x]} (b(r)-\tb(r))\cdot \dee r \\
\begin{split}
&\ge -\int_{p[x^{*},r']} \|b(r)-\tb(r)\|_{2}\|r'(t)\|_{2} \,\dee t \\
&\phantom{\ge~} + \int_{\gamma[x',x]}  (b(r)-\tb(r))\cdot r'(t) \,\dee t 
\end{split} \\
&\ge -(R + \ell^{*})\eps + \int_{\gamma[x',x]} (b(q(t)x)-\tb(q(t)x)) \cdot a x\, \dee t \\
&\ge -(R + \ell^{*})\eps,
\]
where we have used the fact that $r(t) = q(t)x$ for some linear function $q(t)$ with slope $a > 0$. 
Combining the previous two displays, conclude that for all $x \in \mcX$,
$\tpi(x) \le e^{(R+\ell^{*})\eps}\pi(x)$, hence \cref{asm:regularity}(3) holds. 
\eprf

We suspect \cref{prop:checking-integrability} continues to hold even when
$b \ne \grad \log \pi$ and $\tb \ne \grad \log \tpi$. 
Note that condition (1) always holds if $\mcX$ is compact, but also holds for unbounded $\mcX$ as long as
the error in the gradients decays sufficiently quickly as $\|x\|_{2}$ grows large. 
Given an approximate distribution for which $\|b - \tb\|_{2} \le \eps/2$, it is easy to construct
a new distribution that satisfies condition (2):

\bnprop \label{prop:satisfying-condition-2}
Assume that $\tpi$ satisfies $\|b - \tb\|_{2} \le \eps/2$ and let
\[
f_{R}(x) \defined -\frac{\eps x}{2\|x\|_{2}}\left\{(2\|x\|_{2}/R - 1)\ind[R/2 \le \|x\|_{2} < R] + \ind[\|x\|_{2} \ge R]\right\}.
\]
Then the distribution
\[
\tpi_{R}(x) \propto \tpi(x)e^{f_{R}(x)}
\]
satisfies condition (2) of \cref{prop:checking-integrability}.
\enprop

\bprf
Let $\tb_{R} \defined \grad \log \tpi_{R}$. 
First we verify that  $\|b - \tb_{R}\|_{2} \le \eps$.
For $x \in \mcX \cap \ball_{R/2}$, $\tpi_{R}(x) = \tpi(x)$, so $\|b(x) - \tb_{R}(x)\|_{2} \le \eps/2$.
Otherwise $x \in \mcX \setminus \ball_{R/2}$, in which case since $\|f_{R}(x)\| \le \eps/2$
it follows that $\|b(x) - \tb_{R}(x)\|_{2} \le \eps$. 
To verify condition (2), calculate that for $x \in \mcX \setminus \ball_{R}$,
\[
x \cdot (b(x) - \tb_{R}(x))
&= x \cdot \left(b(x) - \tb(x) - \frac{\eps x}{2\|x\|_{2}}\right)
\ge \frac{\eps\|x\|_{2}}{2} - \frac{x \cdot \eps x}{2\|x\|_{2}}
= 0.
\]
\eprf

By taking $R$ very large in \cref{prop:satisfying-condition-2}, we can ensure the integrability condition
holds without having any practical effect on the approximating drift since 
$\tb_{R}(x) = \tb(x)$ for all $x \in \ball_{R/2}$.
Thus, it is safe to view \cref{asm:regularity}(3) as a mild regularity condition. 

\section{Approximation Results for Piecewise Deterministic Markov Processes}
\label{app:pdmps}

In the section we obtain results for a broader class of PDMPs which includes the ZZP a 
special case~\citep{Benaim:2012}.
The class of PDMPs we consider are defined on the space $E \defined \reals^{d} \times \mcB$,
where $\mcB$ is a finite set. 
Let $A \in C^{0}(E, \reals_{+}^{\mcB})$
and let $F \in C^{0}(E, \reals^{d})$ be such that for each 
$\theta \in \mcB$, $F(\cdot, \theta)$ is a smooth vector field
for which the differential equation $\partial_{t}x_{t} = F(x_{t}, \theta)$ with initial condition $x_{0} = x$
has a unique global solution. 
For $\phi \in C(E)$, the standard differential operator $\grad_{x}\phi(x, \theta) \in \reals^{d}$ is
given by $(\grad_{x}\phi(x, \theta))_{i} \defined \D{\phi}{x_{i}}(x, \theta)$ for $i \in [d]$ and
the discrete differential operator $\grad_{\theta}\phi(x, \theta) \in \reals^{\mcB}$
is given by $(\grad_{\theta}\phi(x,\theta))_{\theta'} \defined \phi(x, \theta') - \phi(x, \theta)$. 
The PDMP $(X_{t}, \Theta_{t})_{t\ge0}$ determined by the pair $(F, A)$ has infinitesimal generator
\[
\gen{F,A}\phi = F \cdot \grad_{x} \phi + A \cdot \grad_{\theta} \phi.
\]

We consider the cases when either or both of $A$ and $F$ are approximated 
(in the case of ZZP, only $A$ is approximated while $F$ is exact).
The details of the polynomial contractivity condition depend on which parts of $(F,A)$ 
are approximated. 
We use the same notation for the true and approximating PDMPs with, respectively,
infinitesimal generators $\gen{F,A}$ and $\gen{\tF,\tA}$, as we did for the 
ZZPs in \cref{sec:pdmps}.

\begin{assumption}[\textbf{PDMP error and polynomial contractivity}] \label{asm:pdmp-polynomial-contractivity}
~
\begin{enumerate}
\item There exist $\eps_{F}, \eps_{A} \ge 0$ such that
$\|F - \tF\|_{2} \le \eps_{F}$ and $\|A - \tA\|_{1} \le \eps_{A}$.
\item For each $(x, \theta) \in E$, let $\mu_{x,\theta,t}$ denote the law of the 
PDMP $(X_{x,\theta,t}, \Theta_{x,\theta,t})$ with generator $\gen{F,A}$.
There exist constants $\alpha > 1$ and $\beta > 0$ and a function $B \in C(E \times E, \reals_{+})$
such that for all $x, x' \in \reals^{d}$ and $\theta,\theta' \in \mcB$, 
\[
d_{\mcW}(\mu_{x,\theta,t}, \mu_{x',\theta',t}) &\le B(x, \theta, x', \theta')(t + \beta)^{-\alpha}.
\]
Furthermore, if $\eps_{F} > 0$, then there exists $C_{F} > 0$ such that $B(x, \theta, x', \theta) \le C_{F}\|x - x'\|_{2}$
and if $\eps_{A} > 0$, then there exists $C_{A} > 0$ such that $B(x, \theta, x, \theta') \le C_{A}$.
If $\eps_{F} = 0$ take $C_{F} = 0$ and if $\eps_{A} = 0$ take $C_{A} = 0$. 
\end{enumerate}
\end{assumption}

We also require some regularity conditions similar to those 
for diffusions.

\begin{assumption}[\textbf{PDMP regularity conditions}] \label{asm:pdmp-regularity}
Let $\pi$ and $\tpi$ denote the stationary distributions of the 
PDMPs with, respectively, infinitesimal generators $\gen{F,A}$ and $\gen{\tF,\tA}$. 
\benum
\item The stationary distributions $\pi$ and $\tpi$ exist.
\item The target density satisfies $\int_{E}x^{2}\pi(\dee x, \dee \theta) < \infty$.
\item If a function $\phi \in C(E, \reals)$ is $\pi$-integrable then it is $\tpi$-integrable.
\eenum
\end{assumption}

\bnthm[PDMP error bounds] \label{thm:pdmp-approx-distance}
If \cref{asm:pdmp-polynomial-contractivity,asm:pdmp-regularity} hold, then
\[
d_{\mcW}(\pi, \tpi) \le \frac{C_{F}\eps_{F} + C_{A}\eps_{A}}{(\alpha - 1)\beta^{\alpha - 1}}. \label{eq:pdmp-approx-bound}
\]
\enthm
\bprf[Proof sketch]
For $h \in C_{L}(\reals^{d})$, we need to solve 
\[
h - \pi(h) = \gen{F,A}u.
\]
Similarly to before, the solution is 
\[
u_{h}(x, \theta) \defined \int_{0}^{\infty}(\pi(h) - \EE[h(X_{x,\theta,t})])\,\dee t,
\]
which can be verified as in the diffusion case using Assumptions \ref{asm:pdmp-polynomial-contractivity}(2) and \ref{asm:pdmp-regularity}.
Furthermore, for $x,x' \in \reals^{d}$ and $\theta, \theta' \in \mcB$, by \cref{asm:pdmp-polynomial-contractivity}(2),
\[
|u_{h}(x, \theta) - u_{h}(x', \theta)|
&\le \|h\|_{L}\int_{0}^{\infty}C_{F}\|x - x'\|_{2}(t + \beta)^{-\alpha}\,\dee t \\
&= \frac{C_{F} \|h\|_{L}}{(\alpha - 1)\beta^{\alpha - 1}}\|x - x'\|_{2}
\]
and
\[
|u_{h}(x, \theta) - u_{h}(x, \theta')|
&\le \|h\|_{L}\int_{0}^{\infty}C_{A}(t + \beta)^{-\alpha}\,\dee t 
= \frac{C_{A}\|h\|_{L}}{(\alpha - 1)\beta^{\alpha - 1}}.
\]
We bound $d_{\mcW}(\pi, \tpi)$ as in \cref{thm:stochastic-approx-grad-distance},
but now using the fact that for $u = u_{h}$, $h \in C_{L}(\reals^{d})$, we have
\[
\gen{F,A}u_{h} - \gen{\tF,\tA}u_{h}
&= (F - \tF) \cdot \grad_{x}u_{h} + (A - \tA)\cdot\grad_{\theta}u_{h} \\
&\le \|F - \tF\|_{2}\|\grad_{x}u_{h}\|_{2} + \|A - \tA\|_{1}\|\grad_{\theta}u_{h}\|_{\infty} \\
&\le \frac{C_{F}\eps_{F} + C_{A}\eps_{A}}{(\alpha - 1)\beta^{\alpha - 1}}.
\]
\eprf

\subsection{Hamiltonian Monte Carlo}

We can write an idealized form of Hamiltonian Monte Carlo (HMC) as a PDMP $(X_{t}, P_{t})_{t\ge0}$
 by having the momentum vector $P_{t} \in \reals^{d}$ refresh at a constant rate $\lambda$. 
Let $R_{t}$ be a compound Poisson process with rate $\lambda > 0$ and jump size distribution 
$\distNorm(0, M)$, where $M \in \reals^{d\times d}$ is a positive-definite mass matrix.
That is, if $\Gamma_{t}$ is a homogenous Poisson (counting) process with rate $\lambda$ and 
$J_{i} \dist \distNorm(0, M)$, then
\[
R_{t} \dist \sum_{i=1}^{\Gamma_{t}} J_{i}.
\]
We can then write the HMC dynamics as 
\[
\dee X_{t} &= M^{-1}P_{t}\,\dee t  \\
\dee P_{t} &= \grad \log \pi(X_{T})\, \dee t + \dee R_{t}.
\]
The infinitesimal generator for $(X_{t}, P_{t})_{t\ge0}$ is
\[
\lefteqn{\gen{\lambda,M,\pi}\phi(x, p)} \\
&= M^{-1}p \cdot \grad_{x}\phi(x, p) + \grad \log \pi(x) \cdot \grad_{p}\phi(x,p) + \lambda\left(\int \phi(x, p') \nu_{M}(\dee p') - \phi(x, p)\right),
\]
where $\nu_{M}$ is the density of $\distNorm(0, M)$. 
Let $\mu_{x,p,t}$ denote the law of $(X_{x,p,t}, P_{x,p,t})$ with
generator $\gen{\lambda,M,\pi}$. 
The proof of the following theorem is similar to that for \cref{thm:pdmp-approx-distance}:
\bnthm[HMC error bounds]
Assume that:
\benum
\item $\|\grad\log\pi - \grad\log\tpi\|_{2} \le \eps$.
\item There exist constants $C > 0$ and $0 < \rho < 1$ such that 
\[
d_{\mcW}(\mu_{x,p,t}, \mu_{x,p',t}) \le C\|p - p'\|_{2}\rho^{t}.
\]
\item
The stationary distributions of the PDMPs with, respectively, infinitesimal generators $\gen{\lambda,M,\pi}$ and $\gen{\lambda,M,\tpi}$,
exist (they are, respectively, $\pi \times \mu_{M}$ and $\tpi \times \mu_{M}$).
\item The target density satisfies $\int_{E}x^{2}\pi(\dee x) < \infty$.
\item If a function $\phi \in C(\reals^{d}, \reals)$ is $\pi$-integrable then it is $\tpi$-integrable.
\eenum
Then
\[
d_{\mcW}(\pi, \tpi) \le \frac{C\eps}{\log(1/\rho)}.
\]
\enthm

\section{Analysis of computational--statistical trade-off} 

In this section we prove \cref{thm:ula-accuracy-main}. 
In order to apply results on the approximation accuracy 
of ULA~\citep{Dalalyan:2017,Bubeck:2015,Durmus:2016}, we need
the following property to hold for the exact and approximate drift functions.

\begin{assumption}[Strong log-concavity] \label{asm:strong-log-concavity}
There exists a positive constant $k_{f} > 0$ such that for all $x, x' \in \mcX$,
\[
(f(x) - f(x')) \cdot (x - x') &\le -k_{f} \|x - x'\|_{2}^{2}.
\]
\end{assumption}

We restate the requirements given in \cref{asm:phi-simple}
with some additional notations. 

\begin{assumption}
~
\begin{enumerate}
\item The function $\log \pi_{0} \in C^{3}(\reals^{d}, \reals)$ is $k_{0}$-strongly concave, 
$L_{0} \defined \|\grad \log \pi_{0}\|_{L} < \infty$,
and $\|H[\partial_{j}\log\pi_{0}]\|_{2} \le M_{0} < \infty$ for $j=1,\dots,d$.  
\item There exist constants $k_{\phi}$, $L_{\phi}$, and $M_{\phi}$ such that for $i=1,\dots,N$, 
the function $\phi_{i} \in C^{3}(\reals, \reals)$ is $k_{\phi}$-strongly concave, 
$\|\phi_{i}'\|_{L} \le L_{\phi} < \infty$, and $\|\phi_{i}'''\|_{\infty} \le M_{\phi} < \infty$.
\item The matrix $A_{N} \defined \sum_{i=1}^{N} y_{i}y_{i}^{\top}$ satisfies $\|A_{N}\|_{2} = \Theta(N)$.
\end{enumerate}
\label{asm:phi}
\end{assumption}

Note that under \cref{asm:strong-log-concavity}, there is a unique $x^{\star} \in \reals^{d}$ 
such that $f(x^{\star}) = 0$. 
Our results in this section on based on the following bound on the Wasserstein distance 
between the law of ULA Markov chain and $\pi_{f}$:

\bnthm[{\citep[Theorem 3]{Durmus:2016}, \citep[Corollary 3]{Durmus:2016supp}}] \label{thm:ula-accuracy}
Assume that \ref{asm:strong-log-concavity} holds and the $L_{f} \defined \|f\|_{L} < \infty$. 
Let $\kappa_{f} \defined 2 k_{f}L_{f} / (k_{f} + L_{f})$ and let $\mu_{x,T}$ denote the law of $X_{x,T}'$.
Take $\gamma_{i} = \gamma_{1}i^{-\alpha}$ with $\alpha \in (0,1)$ and set
\[ 
\gamma_{1} = 2(1-\alpha)\kappa_{f}^{-1}(2/T)^{1-\alpha} \log\left(\frac{\kappa_{f} T}{2(1-\alpha)}\right).
\]
If $\gamma_{1} < 1/(k_{f} + L_{f})$, then
\[ 
d_{\mcW}^{2}(\mu_{x,T}, \pi_{f})
\le 16(1 - \alpha)L_{f}^{2}\kappa_{f}^{-3} d T^{-1}\log\left(\frac{\kappa_{f} T}{2(1-\alpha)}\right).
\]
\enthm

For simplicity we fix $\alpha = 1/2$, though the same results hold for all $\alpha \in (0,1)$, 
just with different constants. 
Take $\{\gamma_{i}\}_{i=1}^{\infty}$ as defined in \cref{thm:ula-accuracy}. 
Let $x^{\star} \defined \argmax_{x} \mcL(x)$ and let $S_{k} \defined \sum_{i=1}^{N} \|y_{i}\|_{2}^{k}$. 
The drift for this model is
\[
b(x) \defined \grad \mcL(x) =  \grad \log \pi_{0}(x) + \sum_{i=1}^{N} \phi_{i}'(x \cdot y_{i}) y_{i}.
\]
By Taylor's theorem, the $j$-th component of $b(x)$ can be rewritten as 
\[
\begin{split}
b_{j}(x) 
&= \partial_{j} \log \pi_{0}(x^{\star}) + \grad \partial_{j}\log \pi_{0}(x^{\star})\cdot(x - x^{\star}) + R(\partial_{j} \log \pi_{0}, x) \\
&\phantom{=~} + \sum_{i=1}^{N} \phi_{i}'(x^{\star} \cdot y_{i}) y_{ij} + \phi_{i}''(x^{\star} \cdot y_{i})y_{ij}y_{i}\cdot(x - x^{\star}) + R(\phi_{i}'(\cdot \cdot y_{i})y_{ij}, x) 
\end{split} \\
\begin{split}
&= \grad \partial_{j}\log \pi_{0}(x^{\star})\cdot(x - x^{\star}) + R(\partial_{j} \log \pi_{0}, x) \\
&\phantom{=~}  + \sum_{i=1}^{N} \phi_{i}''(x^{\star} \cdot y_{i})y_{ij}y_{i}\cdot(x - x^{\star}) + R(\phi_{i}'(\cdot \cdot y_{i})y_{ij}, x),
\label{eq:b-taylor-expansion}
\end{split} 
\]
where 
\[
R(f, x) \defined \|x - x^{\star}\|_{2}^{2}\int_{0}^{1}(1-t)Hf(x^{\star} + t(x - x^{\star}))\,\dee t.
\]
Hence we can approximate the drift with a first-order Taylor expansion around $x^{\star}$:
\[
\tb(x) \defined  (H\log\pi_{0})(x^{\star})(x - x^{\star}) + \sum_{i=1}^{N} \phi_{i}''(x^{\star} \cdot y_{i}) y_{i}y_{i}^{\top}(x - x^{\star}).
\]
Observe that \cref{asm:strong-log-concavity} is satisfied for $f = b$ and $f = \tb$ with 
$k_{f} = k_{N} \defined k_{0} + k_{\phi}\|A_{N}\|_{2}$.
Furthermore, \cref{asm:regularity} is satisfied with $\|\tb\|_{L} \le L_{N} \defined L_{0} + L_{\phi}S_{2}$
and $\|b\|_{L} \le L_{N}$ as well since
\[
\|\phi_{i}'(x_{1} \cdot y_{i})y_{i} - \phi_{i}'(x_{2} \cdot y_{i})y_{i}\|_{2}
&\le |\phi_{i}'(x_{1} \cdot y_{i}) - \phi_{i}'(x_{2} \cdot y_{i})| \|y_{i}\|_{2} \\
&\le L_{\phi}|x_{1} \cdot y_{i} - x_{2} \cdot y_{i}| \|y_{i}\|_{2} \\
&\le L_{\phi} \|y_{i}\|_{2}^{2} \|x_{1} - x_{2}\|_{2}.
\]
Thus, $b$ and $\tb$ satisfy the same regularity conditions.%

We next show that they cannot deviate too much from each other.
Using \cref{eq:b-taylor-expansion} and regularity assumptions we have
\[
\|b(x) - \tb(x)\|_{2}^{2}
&= \sum_{j=1}^{d}\left(R(\partial_{j} \log \pi_{0}, x) + \sum_{i=1}^{N} R(\phi_{i}'(\cdot \cdot y_{i})y_{ij}, x)\right)^{2} \\
&\le \|x - x^{\star}\|_{2}^{4}\sum_{j=1}^{d}\left(M_{0} + \sum_{i=1}^{N} M_{\phi}\|y_{i}\|_{2}^{2}y_{ij} \right)^{2} \\
&\le d\|x - x^{\star}\|_{2}^{4}\left(M_{0} + M_{\phi}\sum_{i=1}^{N} \|y_{i}\|_{2}^{3}\right)^{2}. 
\] 
It follows from \citep[Theorem 1(ii)]{Durmus:2016} that
\[
\tpi(\|b - \tb\|_{2}) \le d^{3/2}M_{N}k_{N}^{-1},
\]
where $M_{N} \defined M_{0} + M_{\phi}S_{3}$.

Putting these results together with \cref{thm:ula-accuracy,thm:approx-grad-distance}
and applying the triangle inequality, we conclude that
\[
d_{\mcW}^{2}(\mu^{\star}_{T}, \pi) 
&\le \frac{(k_{N} + L_{N})^{3}d}{k_{N}^{3}L_{N}} \frac{\log\left(\frac{2 k_{N}L_{N}}{k_{N} + L_{N}}T\right)}{T} \\
d_{\mcW}^{2}(\tmu^{\star}_{\tT}, \pi) 
&\le \frac{2(k_{N} + L_{N})^{3}d}{k_{N}^{3}L_{N}} \frac{\log\left(\frac{2 k_{N}L_{N}}{k_{N} + L_{N}}\tT\right)}{\tT} + \frac{2d^{3}M_{N}^{2}}{k_{N}^{4}}.
\]
In order to compare the bounds we must make the computational budgets of the two algorithms equal.
Recall that we measure computational cost by the number of $d$-dimensional inner products performed, so ULA with $b$ 
costs $TN$ and ULA with $\tb$ costs $(\tT + N)d$. 
Equating the two yields $\tT = N(T/d - 1)$, so we must assume that $T > d$. 
For the purposes of asymptotic analysis, assume also that $S_{i}/N$ is bounded from above and bounded
away from zero. 
Under these assumptions, in the case of $k_{\phi} > 0$, we conclude that
\[
d_{\mcW}^{2}(\mu^{\star}_{T}, \pi) = \tO\left(\frac{d}{TN}\right)
\qquad \text{and} \qquad
d_{\mcW}^{2}(\tmu^{\star}_{\tT}, \pi) = \tO\left(\frac{d^{2}}{N^{2}T} + \frac{d^{3}}{N^{2}}\right),
\]
establishing the result of \cref{thm:ula-accuracy-main}. For large $N$, the approximate ULA with $\tb$ is more accurate.

\bibliographystyle{abbrvnat}
\bibliography{library}

\end{document}